\documentclass[10pt]{amsart}

\usepackage{geometry}
\usepackage{amsmath}
\usepackage{amssymb}
\usepackage[foot]{amsaddr}
\usepackage{enumitem}
\usepackage{cases}

\calclayout

\newcommand{\R}{\mathbb{R}}

\newcommand{\PP}{\mathbb{P}}
\newcommand{\FF}{\mathcal{F}}
\newcommand{\E}{\mathbb{E}}

\def\qed{\hfill \hbox{$\Box$} \smallskip}

\newtheorem{theorem}{Theorem}[section]
\newtheorem{lemma}[theorem]{Lemma}
\newtheorem{proposition}[theorem]{Proposition}
\newtheorem{assumption}[theorem]{Assumption}
\newtheorem{corollary}[theorem]{Corollary}

\newtheorem{definition}[theorem]{Definition}

\newtheorem{remark}[theorem]{Remark}

\numberwithin{equation}{section}

\begin{document}  

\title[Singular Perturbation of Stochastic Differential Games]{Singular Perturbation of Zero-Sum Linear-Quadratic Stochastic Differential Games}

\author{Beniamin Goldys}
\author{James Yang}
\author{Zhou Zhou}
\address{School of Mathematics and Statistics, University of Sydney, Australia.}

\date{\today}
\maketitle
\begin{abstract}

We investigate a class of zero-sum linear-quadratic stochastic differential games on a finite time horizon governed by multiscale state equations. The multiscale nature of the problem can be leveraged
to reformulate the associated generalised Riccati equation as a deterministic singular perturbation problem. In doing so, we show that, for small enough $\epsilon$, the existence of solution to the associated generalised Riccati equation is guaranteed by the existence of a solution to a decoupled pair of differential and algebraic Riccati equations with a reduced order of dimensionality. Furthermore, we are able to formulate a pair of asymptotic estimates to the value function of the game problem by constructing an approximate feedback strategy and observing the limiting value function.

\end{abstract}

\section{Introduction}

In the paper, we study a class of two-player stochastic differential games on a finite time interval $[0,T]$. The dynamical system is described by a slow process $X_1$ and a fast process $X^\epsilon_2$ that satisfy the following linear stochastic differential equations
\begin{align*}
\begin{cases}
dX_1(t) = \left[A_{11} X_1(t) + A_{12} X^\epsilon_2(t)  + B_{11} u_1(t) + B_{12} u_2(t) \right]dt + \sigma_1 dW_1(t),\\
dX^\epsilon_1(t) = \frac{1}{\epsilon} \left[ A_{21} X(t) + A_{22} X^\epsilon_2(t) + B_{21} u_1(t) + B_{22} u_2(t) \right] dt + \frac{\sigma_2}{\sqrt{\epsilon}} dW_2(t),\\
X_1(0) = x_1,\quad X^\epsilon_2(0) = x_2.
\end{cases}
\end{align*}
Here $\epsilon$ is a small positive parameter representing the ratio between the evolutionary speeds of the slow and fast processes. The objective of the two-player game adheres to a zero-sum formulation. That is, Player 1 attempts to maximise a quadratic objective function through a strategic decision $u_1$, whilst Player 2 attempts to minimise the same quadratic objective function through a strategic decision $u_2$. A classical approach to obtaining a feedback saddle point or equilibrium is to consider the associated generalised Riccati equation. Our interests centre on the solvability of the generalised Riccati equation and obtaining asymptotic estimates for the value function of the game when $\epsilon$ is small. Due to the degenerate nature of the differential term when $\epsilon$ is formally set to be $0$, care needs to be applied when studying the convergence problem. Problems with this asymptotic property are often referred to as \textit{singular perturbation problems}. Suitably, we refer to the game problem described in this paper as a \textit{singularly perturbed zero-sum linear-quadratic stochastic differential game}, or Problem (SLQG) for short.

In 1965, the study of linear-quadratic differential games in the zero-sum framework was initiated by Ho, Bryson and Baron \cite{ho1965differential} in the context of pursuit-evasion strategies. As opposed to what is observed in optimal control theory, open-loop and feedback saddle points are not necessarily the same nor do they imply the other exists. This important distinction was first brought to light by Schmitendorf \cite{schmitendorf1970existence} and was later shown by Sun and Yong \cite{sun2014linear} to be equivalent as long as the saddle points exist and the objective function possesses a certain convexity-concavity assumption. In 1979, Bernhard \cite{bernhard1979linear} examined in great depth the case of feedback strategies in deterministic differential games. In the stochastic framework, Sun and Yong \cite{sun2014linear} demonstrated that the existence of a feedback saddle point is equivalent to the solvability of a generalised Riccati equation. Due to the indefinite nature of the nonlinear terms, obtaining analytically checkable conditions for the solvability of the generalised Riccati equation has proved to be a challenge in and among itself. Yong \cite{yong1999linear} showed that the existence a solution to the generalised Riccati equation is equivalent to the invertibility of a submatrix of a matrix exponential. On the other hand, McAsey and Mou \cite{mcasey2006generalized} derived a comparison theorem and an existence result based on the existence of upper and lower solutions. Under the framework of non-anticipating strategies, Yu \cite{yu2015optimal} stated commutative and definiteness conditions that prove to be sufficient for the solvability of the generalised Riccati equation. Other particular cases, such as the one dimensional case, and numerical examples are discussed in the book by Yong \cite{jiongmin2014differential}. 

In 1977, the concept of singular perturbation in zero-sum linear-quadratic differential games was first considered by Gardner Jr. \cite{gardner1977zero}. Under a finite time horizon and a deterministic setting, the author constructs a composite pair of strategies from the saddle points of slow and fast sub-games to obtain an approximate value function to the original game. In the infinite time horizon and non-zero sum counterpart, a composite strategy approach was also adopted by Gardner Jr. and Cruz Jr. \cite{gardner1978well} and Khalil and Kokotovic \cite{khalil1979feedback} to produce similar asymptotic estimates. The non-linear case has been studied, via a dynamic programming approach, in the deterministic case by Gaitsgory \cite{gaitsgory1996limit} and Subbotina \cite{subbotina1996asymptotic, subbotina1999asymptotic, subbotina2001asymptotics} and in the stochastic case by Alvarez and Bardi \cite{alvarez2010ergodicity}. In this case, a singularly perturbed Hamilton-Jacobi-Bellman equation is formulated and analysed. Differential games are a natural extension of optimal control problems, which suitably can be refered to as a single player game. That being said, applications of singular perturbation in optimal control theory are plentiful, whether the motivation is to describe a multiscale optimal control problem or to reduce the order of complexity. For an extensive survey of applications to areas such as aerospace engineering, biology and chemistry, see \cite{naidu1988singular, naidu2002singular} and for stochastic filtering in finance, see \cite{fouque2015filtering, fouque2017perturbation}.

In this paper, we make two contributions, both of which flesh out a framework for zero-sum linear-quadratic stochastic differential games with slow-fast state processes. First, we show that for any finite time horizon, if $\epsilon$ is small enough then a sufficient condition for the solvability of the generalised Riccati equation is the solvability of a decoupled pair of differential and algebraic Riccati equations with a lower order of dimensionality. More precisely, we reduce the original problem of order $(n_1+n_2)^2$ to a pair of problems of orders $n_1^2$ and $n_2^2$. Using this, the feedback saddle point and value function of Problem (SLQG) can be shown to exist. Secondly, we produce an asymptotic estimate to the value function of Problem (SLQG) of order $O(\epsilon)$ by constructing an approximate feedback saddle point using feedback operators when $\epsilon$ is formally set to be zero. In addition, we show that the value function of Problem (SLQG) converges to a limiting value function.

We approach Problem (SLQG) by characterising the feedback saddle point and value function in terms of the solution to a generalised Riccati equation, see Sun and Yong \cite{sun2014linear}. To eliminate the singular terms in the generalised Riccati equation, we suppose that the solution adopts a first-order representation \cite{dragan2012linear, sannuti1969near}. As a result, we obtain a system of ODEs, which we refer to as the \textit{full system}. At this point, it can be observed that the full system represents a classical and deterministic singular perturbation problem. To this end, we consider the \textit{reduced system} of differential and algebraic equations when $\epsilon$ is formally set to zero. The solvability of the reduced system is shown to be guaranteed by the solvability of a decoupled pair of differential and algebraic Riccati equations of a lower dimension \--- which we shall refer to as the \textit{reduced differential-algebraic Riccati equations}. By formulating and proving a version of the Tikhonov Theorem, we establish that for sufficiently small $\epsilon$, the solvability of the reduced differential-algebraic Riccati equations implies the existence of a unique solution to the full system and that its limit exists. This the key result of this paper and from this, we can deduce our main results outlined in the previous paragraph. 

To the best of our knowledge, singular perturbations in linear-quadratic differential games on a finite time horizon has not been considered in the stochastic framework. In the deterministic case, Gardner Jr. \cite{gardner1977zero} uses a composite control approach and assume the solution of the generalised Riccati equation to exist and possess analytic properties, which we do not require. On the other hand, the non-linear counterparts \cite{alvarez2010ergodicity, gaitsgory1996limit, subbotina1996asymptotic, subbotina1999asymptotic, subbotina2001asymptotics} require that the objective function be uniformly Lipscthiz with respect to the state variables, which excludes the quadratic case. We remark that while we assume the coefficient matrices are constant, there is no mathematical difficulty in extending the methods and results to bounded and deterministic time-dependent coefficients. Lastly, we assume that the terminal objective function is zero, else otherwise the full system of ODEs will have unbounded terminal conditions. Problems of this class have been studied by Glizer \cite{glizer1997asymptotic, glizer2000asymptotic}.

The paper is organised as follow: In Section \ref{section: formulation}, we formulate the singularly perturbed zero-sum linear-quadratic stochastic differential game and introduce the associated generalised Riccati equation. In Section \ref{section: reduced}, we make the connection between the reduced system and the reduced differential-algebraic Riccati equations. Section \ref{section: asymptotic} sits at the heart of this paper, establishing the existence, uniqueness and convergence of the solution of the generalised Riccati equation via a version of the Tikhonov theorem.  Section \ref{section: approxcontrol} is devoted to obtaining asymptotic estimates for the value function of the game. In Section \ref{section: onedimension}, we summarise our results for the one dimensional case. Finally, in Section \ref{section: proof}, we give the proof for the Tikhonov theorem used in this paper.

\section{Mathematical formulation}\label{section: formulation}

\subsection{Notation}

Given a real and separable Hilbert space $E$, the inner product of its elements is denoted by $\langle \cdot, \cdot \rangle_E$ and the associated norm as $|\cdot |_E$. If $G = E\times F$ is the Cartesian product of the Hilbert spaces $E$ and $F$ then $G$ endowed with the inner product $\langle \cdot,\cdot \rangle_G = \langle \cdot,\cdot \rangle_E + \langle \cdot,\cdot \rangle_F$ is also a Hilbert space. For the most part, when there is no confusion, we will drop the subscript in the inner product and norm. For matrices in the space $\R^{n\times n}$, we will write $\mathbb{S}^n$ as the space of symmetric matrices, $\mathbb{S}^n_+$ as the space of non-negative (or positive semi-definite) symmetric matrices and $\mathbb{S}^n_{++}$ as the space of positive (definite) symmetric matrices. For a matrix $M$, we shall use the notation $M^*$ to denote the transpose of $M$ and $M^{-*}$ to denote the inverse of the transpose.

Let $T > 0$ be a finite time horizon and let  $W_1 := (W_1(t))_{0\leq t\leq T}$ and $W_2 := (W_2(t))_{0\leq t\leq T}$ be two independent $m_1$ and $m_2$-dimensional Brownian motions, respectively, defined on a complete filtered probability space $(\Omega, \mathcal{F} ,\mathbb{P})$, where $(F_t)$ is the natural filtration of $W_1$ and $W_2$ augmented by the $\PP$-null sets of $\FF$. We define the following spaces of processes with respect to the Hilbert space $E$:
\begin{itemize}
\item $C([0,T]; E)$ is the space of continuous mappings $F: [0,T] \rightarrow E$ equipped with the norm
\[
\| F\|_{C([0,T]; E)} := \sup_{t\in [0,T]} |F(t)|_{E}.
\] 
\item $L^2_\mathcal{F}(\Omega \times [0,T]; E)$ is the space of equivalence classes of processes $F \in L^2(\Omega \times [0,T]; E)$ admitting a predictable version and equipped with the norm
\[
\| F\|_{L^2_\mathcal{F}(\Omega \times [0,T]; E)} := \left( \E\int_0^T |F(t)|_{E}^2 dt \right)^{1/2}.
\]
For short, we will write $\mathbb{A}_T^n := L^2_\mathcal{F}(\Omega \times [0,T]; \R^n)$.
\item $L^2_\mathcal{F}(\Omega; C([0,T];E))$ is the space of predictable processes $F: \Omega \times [0,T] \rightarrow E$ with continuous paths in $E$ equipped with the norm
\[
\| F\|_{L^2_\mathcal{F}(\Omega; C([0,T];E))} := \left( \E \sup_{t\in [0,T]} |F(t)|_{E}^2  \right)^{1/2}.
\]
For short, we will write $\mathbb{H}_T^n := L^2_\mathcal{F}(\Omega; C([0,T];\R^n))$.
\end{itemize}

\subsection{Problem formulation}

Let $0 < \epsilon \leq 1$ denote the perturbation parameter. For arbitrary $x_1 \in \R^{n_1}$ and $x_2 \in \R^{n_2}$, the slow state $X_1 := (X_1(t))_{0\leq t\leq T}$ and fast state $X^\epsilon_2 := (X^\epsilon_2(t))_{0\leq t\leq T}$ take values in $\R^{n_1}$ and $\R^{n_2}$ respectively, and evolve according to the following system of linear stochastic differential equations
\begin{equation}\label{state}
\begin{split}
\begin{cases}
dX_1(t) = \left[A_{11} X_1(t) + A_{12} X^\epsilon_2(t)  + B_{11} u_1(t) + B_{12} u_2(t) \right]dt + \sigma_1 dW_1(t),\\
dX^\epsilon_2(t) = \frac{1}{\epsilon} \left[ A_{21} X_1(t) + A_{22} X^\epsilon_2(t) + B_{21} u_1(t) + B_{22} u_2(t) \right] dt + \frac{\sigma_2}{\sqrt{\epsilon}} dW_2(t),\\
X_1(0) = x_1,\quad X^\epsilon_2(0) = x_2.
\end{cases}
\end{split}
\end{equation}
Here the coefficients $A_{ij}, B_{ij}, \sigma_i, i,j = 1,2,$ are constant matrices of appropriate dimensions. Player 1 (resp. Player 2) influences the state equation with respect to the strategy $u_1 := (u_1(t))_{t\in [0,T]}$ (resp. $u_2:=(u_2(t))_{t\in [0,T]}$) taking values in $\R^{k_1}$ (resp. $\R^{k_2}$). The objective of Player 1 (resp. Player 2) is to maximise (resp. minimise) the following quadratic objective functional with respect to an admissible strategy $u_1 \in \mathbb{A}_T^{k_1} := L^2_\mathcal{F}(\Omega\times [0,T]; \R^{k_1})$ (resp. $u_2 \in \mathbb{A}_T^{k_2}:= L^2_\mathcal{F}(\Omega\times [0,T]; \R^{k_2})$)
\begin{equation}\label{cost}
J^\epsilon(x_1,x_2; u_1, u_2) = \frac{1}{2} \E \int_0^T \left[ \left\langle Q_1 X_1(t), X_1(t) \right\rangle + \left\langle Q_2 X^\epsilon_2 (t),X^\epsilon_2 (t) \right\rangle   -  |u_1(t)|^2 + |u_2(t)|^2\right]dt.
\end{equation}
where $Q_1$ and $Q_2$ are constant symmetric matrices of appropriate dimensions. 

Let $n = n_1 + n_2$. Let us introduce the following compact notation:
\begin{align*}
\begin{cases}
X^\epsilon(t) = \begin{pmatrix} X_1(t)\\ X^{\epsilon}_2(t) \end{pmatrix}, \quad u(t) = \begin{pmatrix} u_1(t) \\ u_2(t) \end{pmatrix}, \quad W(t) = \begin{pmatrix} W_1(t) \\ W_2(t) \end{pmatrix}, \quad Q = \begin{pmatrix}
Q_1 & 0\\
0 & Q_2
\end{pmatrix}, \\
A^\epsilon = \begin{pmatrix} A_{11} & A_{12} \\ \frac{1}{\epsilon} A_{21} & \frac{1}{\epsilon}A_{22}  \end{pmatrix},\quad 
B^\epsilon = \begin{pmatrix} B_{11} & B_{12} \\ \frac{1}{\epsilon}B_{21} & \frac{1}{\epsilon} B_{22} \end{pmatrix},\quad \sigma^\epsilon = \begin{pmatrix} \sigma_1 & 0 \\ 0 & \frac{1}{\sqrt{\epsilon}}\sigma_2 \end{pmatrix}, \quad R = \begin{pmatrix}
-I & 0\\
0 & I
\end{pmatrix}.
\end{cases}
\end{align*} 
As a result, we can rewrite the state equation \eqref{state} and the objective function \eqref{cost} as 
\begin{equation}\label{statecompact}
\begin{split}
\begin{cases}
dX^\epsilon(t) = \left[ A^\epsilon X^\epsilon(t) + B^\epsilon u(t) \right] dt + \sigma^\epsilon dW(t),\\
X^\epsilon(0) = x := \begin{pmatrix}
x_1\\x_2
\end{pmatrix},
\end{cases}
\end{split}
\end{equation}
and 
\begin{equation}\label{costcompact}
J^\epsilon(x; u) = \frac{1}{2} \E \int_0^T  \left[ \langle Q X^\epsilon(t) , X^\epsilon(t) \rangle + \langle R u(t), u(t) \rangle \right] dt
\end{equation}
respectively.
For fixed $0<\epsilon\leq 1$, if $(u_1,u_2) \in \mathbb{A}^{k_1}_T\times \mathbb{A}_T^{k_2}$ then the linear state equation \eqref{statecompact} admits a unique solution $X^\epsilon \in \mathbb{H}^n_T := L^2_\FF(\Omega; C([0,T];\R^n))$. Consequently, the coupled state equations \eqref{state} also admit a unique solution $(X_1,X^\epsilon_2) \in \mathbb{H}^{n_1}_T \times \mathbb{H}^{n_2}_T$. Moreover, the objective functions \eqref{cost} and \eqref{costcompact} are well-defined for pairs of admissible strategies $(u_1,u_2) \in \mathbb{A}_T^{k_1} \times \mathbb{A}_T^{k_2}$. We will refer to the game problem characterised by \eqref{state} and \eqref{cost} or \eqref{statecompact} and \eqref{costcompact} as the \textit{singularly perturbed zero-sum linear-quadratic stochastic differential game} or Problem (SLQG) for short. 

In this paper, the notion of a solution is defined in terms of the Nash equilibrium or, more precisely in the context of zero-sum games, a saddle point. In this situation, neither player can non-cooperatively improve their position by choosing an alternative strategy. Furthermore, we shall confine the players to the class of feedback strategies. That is, players make decisions based on an instantaneous and perfect knowledge of the current state: For fixed $\epsilon \in (0,1]$
\begin{equation}\label{closedloopcontrol}
u_i^\epsilon(t) = \phi_i(t,X_1(t),X^\epsilon_2(t), \epsilon),\quad i = 1,2,\ \forall t\in [0,T],
\end{equation}
where $\phi_i(t,X_1,X_2,\epsilon)$ is some deterministic and measurable function that is Lipschitz continuous in state variables $(X_1,X_2)$, uniformly in $t$. Let us denote by $\mathbb{A}_{T,f}^{k_i}$ as the set of strategies, which are of the feedback form \eqref{closedloopcontrol}.
\begin{definition}
Fix $0<\epsilon\leq 1$. A pair of strategies $(\widehat{u}_1^\epsilon, \widehat{u}_2^\epsilon)$ defined as
\begin{align*}
\begin{cases}
\widehat{u}_1^\epsilon(t) = \widehat{\phi}_1(t,\widehat{X}^\epsilon_1(t),\widehat{X}^\epsilon_2(t),\epsilon), \quad t\in [0,T],\\
\widehat{u}_2^\epsilon(t) = \widehat{\phi}_2(t,\widehat{X}^\epsilon_1(t),\widehat{X}^\epsilon_2(t),\epsilon), \quad t\in [0,T],
\end{cases}
\end{align*}
where $\widehat{X}^\epsilon_1(t) := X_1(t;x_1, \widehat{u}_1^\epsilon, \widehat{u}_2^\epsilon)$ and $\widehat{X}^\epsilon_2(t) := X^\epsilon_2(t;x_2, \widehat{u}_1^\epsilon, \widehat{u}_2^\epsilon)$, is a feedback saddle point of Problem (SLQG) if 
\begin{equation}\label{saddledefinition}
J^\epsilon(x_1,x_2; u_1^\epsilon, \widehat{u}_2^\epsilon) \leq J^\epsilon(x_1,x_2; \widehat{u}_1^\epsilon, \widehat{u}_2^\epsilon) \leq J^\epsilon(x_1,x_2; \widehat{u}_1^\epsilon,u_2^\epsilon), \ \forall u_1^\epsilon \in \mathbb{A}^{k_1}_{T},\forall u_2^\epsilon \in \mathbb{A}^{k_2}_{T}.
\end{equation}
\end{definition}
It should be pointed out that the pair $(\widehat{u}_1^\epsilon,u_2^\epsilon)$ in $J^\epsilon(x_1,x_2; \widehat{u}_1^\epsilon,u_2^\epsilon)$ refers to 
\begin{align*}
\begin{cases}
\widehat{u}_1^\epsilon(t) = \widehat{\phi}_1(t,X^\epsilon_1(t),X^\epsilon_2(t),\epsilon), \quad t\in [0,T],\\
u_2^\epsilon(t) = \phi_2(t,X^\epsilon_1(t),X^\epsilon_2(t),\epsilon), \quad t\in [0,T],
\end{cases}
\end{align*}
where $X^\epsilon_1(t) := X_1(t;x_1 \widehat{u}_1^\epsilon, u_2^\epsilon)$ and $X^\epsilon_2(t) := X^\epsilon_2(t;x_2, \widehat{u}_1^\epsilon, u_2^\epsilon)$. In other words, when Player 2 changes his strategy from $\widehat{u}_2^\epsilon$ to $u_2^\epsilon$, Player 1's strategy $\widehat{u}_1^\epsilon$ is impacted due to his or her dependence on the state variables. Similar can be said for the pair $(u_1^\epsilon, \widehat{u}_2^\epsilon)$ in $J^\epsilon(x_1,x_2; u_1^\epsilon, \widehat{u}_2^\epsilon)$. Next, we define the value of Problem (SLQG).
\begin{definition}
The \textit{upper value} $V_+^\epsilon(x_1,x_2)$ and \textit{lower value} $V_-^\epsilon(x_1,x_2)$ are defined as 
\begin{equation}
V_+^\epsilon(x_1,x_2) = \inf_{u_2 \in \mathbb{A}^{k_2}_{T}} \sup_{u_1 \in \mathbb{A}^{k_1}_{T}} J^\epsilon(x_1,x_2;u_1,u_2)
\end{equation}
\begin{equation}
V_-^\epsilon(x_1,x_2) = \sup_{u_1 \in \mathbb{A}^{k_1}_{T}} \inf_{u_2 \in \mathbb{A}^{k_2}_{T}}  J^\epsilon(x_1,x_2;u_1,u_2)
\end{equation}
Clearly, $V_-^\epsilon(x_1,x_2) \leq V_+^\epsilon(x_1,x_2)$. If $V_-^\epsilon(x_1,x_2) = V_+^\epsilon(x_1,x_2) := V^\epsilon(x_1,x_2)$ then we say that Problem (SLQG) admits a \textit{value}.
\end{definition}

A classical approach to feedback differential game problems is by characterising the saddle point and value function in terms of the solution of a generalised Riccati equation
\begin{equation}\label{riccati}
\begin{split}
\begin{cases}
\dot{P}^\epsilon + (A^\epsilon)^* P^\epsilon + P^\epsilon A^\epsilon - P^\epsilon B^\epsilon R^{-1} (B^{\epsilon})^* P^\epsilon  + Q= 0,\\
P^\epsilon(T) = 0.
\end{cases}
\end{split}
\end{equation}
For convenience, we write the non-linear term as
\[
-B^\epsilon R^{-1} (B^\epsilon)^* = \begin{pmatrix}
\Delta_1 & \frac{\Delta}{\epsilon}\\
\frac{\Delta^*}{\epsilon} & \frac{\Delta_2}{\epsilon^2}
\end{pmatrix},
\] 
where 
\begin{align*}
\begin{cases}
\Delta_1 = B_{11} B_{11}^* - B_{12} B_{12}^*, \\
\Delta = B_{11} B_{21}^* - B_{12} B_{22}^*,\\
\Delta_2 = B_{21} B_{21}^* - B_{22} B_{22}^*.
\end{cases}
\end{align*}
The following theorem demonstrates that the existence and uniqueness of a solution to the Riccati equation implies the existence and uniqueness of a feedback saddle point and value.

\begin{theorem}\label{theorem: solution}
Fix $0 < \epsilon\leq 1$. Suppose the generalised Riccati equation admits a unique solution $P^\epsilon \in C([0,T];\mathbb{S}^{n})$. We write $P^\epsilon$ as the first order representation
\[
\begin{pmatrix}
P^\epsilon_{11} & \epsilon P^\epsilon_{12}\\
\epsilon (P^\epsilon_{12})^* & \epsilon P^\epsilon_{22}
\end{pmatrix}
\]
where $P^\epsilon_{11} \in C([0,T];\mathbb{S}^{n_1}), P^\epsilon_{12} \in C([0,T];\R^{n_1\times n_2})$ and $P^\epsilon_{22} \in C([0,T];\mathbb{S}^{n_2})$.
Then Problem (SLQG) admits a unique feedback saddle point defined on the interval $[0,T]$
\begin{equation}\label{saddle}
\begin{split}
\begin{cases}
\widehat{u}^\epsilon_1(t) &= \widehat{\phi}_1(t,\widehat{X}^\epsilon_1(t), \widehat{X}^\epsilon_2(t),\epsilon) \\
&:=  \left[B_{11}^* P^\epsilon_{11}(t) + B_{21}^* (P^\epsilon_{12}(t))^*\right] \widehat{X}^\epsilon_1(t) + \left[\epsilon B_{11}^* P^\epsilon_{12}(t) + B_{21}^* P^\epsilon_{22}(t)\right] \widehat{X}^\epsilon_2(t),\\
\widehat{u}^\epsilon_2(t) &=\widehat{\phi}_2(t,\widehat{X}^\epsilon_1(t), \widehat{X}^\epsilon_2(t),\epsilon)\\
 &:=  - \left[B_{12}^* P^\epsilon_{11}(t) +  B_{22}^* (P^\epsilon_{12}(t))^*\right] \widehat{X}^\epsilon_1(t) - \left[ \epsilon B_{12}^* P^\epsilon_{12}(t) + B_{22}^* P^\epsilon_{22}(t)\right] \widehat{X}^\epsilon_2(t),
\end{cases}
\end{split}
\end{equation}
where $(\widehat{X}^\epsilon_1,\widehat{X}^\epsilon_2)$ is the solution to the state equation
\begin{equation}\label{optimstate}
\begin{split}
\begin{cases}
d\widehat{X}^\epsilon_1(t) = \left[A_{11} \widehat{X}^\epsilon_1(t) + A_{12} \widehat{X}^\epsilon_2(t)  + B_{11} \widehat{\phi}_1(t,\widehat{X}^\epsilon_1(t), \widehat{X}^\epsilon_2(t),\epsilon) + B_{12} \widehat{\phi}_2(t,\widehat{X}^\epsilon_1(t), \widehat{X}^\epsilon_2(t),\epsilon) \right]dt \\
\qquad \qquad + \sigma_1 dW_1(t),\\
d\widehat{X}^\epsilon_2(t) = \frac{1}{\epsilon} \left[ A_{21} \widehat{X}^\epsilon_1(t) + A_{22} \widehat{X}^\epsilon_2(t) + B_{21} \widehat{\phi}_1(t,\widehat{X}^\epsilon_1(t), \widehat{X}^\epsilon_2(t),\epsilon) + B_{22} \widehat{\phi}_2(t,\widehat{X}^\epsilon_1(t), \widehat{X}^\epsilon_2(t),\epsilon) \right] dt \\
\qquad \qquad + \frac{\sigma_2}{\sqrt{\epsilon}} dW_2(t),\\
\widehat{X}^\epsilon_1(0) = x_1,\quad \widehat{X}^\epsilon_2(0) = x_2.
\end{cases}
\end{split}
\end{equation}
and has the value
\begin{equation}
\begin{split}
V^\epsilon(x_1,x_2) &= J^\epsilon(x_1,x_2;\widehat{u}_1,\widehat{u}_2) \\
&=\frac{1}{2}\left[\left\langle P^\epsilon_{11}(0) x_1 ,x_1\right\rangle + 2 \epsilon\left\langle x_1 , P^\epsilon_{12}(0) x_2\right\rangle + \epsilon \left\langle P^\epsilon_{22}(0) x_2 ,x_2\right\rangle\right]  \\
&\quad + \frac{1}{2} \int_0^T \left[\left\langle P^\epsilon_{11}(t) \sigma_1, \sigma_1 \right\rangle + \left\langle P^\epsilon_{22}(t) \sigma_2,\sigma_2 \right\rangle\right]  dt.
\end{split}
\end{equation}
\end{theorem}
\noindent \textit{Proof.}  Similar to the completion of squares argument in Theorem 6.1 of \cite{yong1999stochastic}, we can apply Ito's lemma to $\left\langle P^\epsilon(t) X^\epsilon(t), X^\epsilon(t) \right\rangle$ to show that
\begin{equation}\label{temp0}
\begin{split}
&J^\epsilon (x_1,x_2; u_1,u_2)\\
&= - \frac{1}{2}\E\int_0^T\left[ |u_1(t) - \widehat{F}^\epsilon_{11}(t) X_1(t) - \widehat{F}^\epsilon_{12}(t) X^\epsilon_2(t) |^2 - |u_2(t) - \widehat{F}^\epsilon_{21}(t) X_1(t) - \widehat{F}^\epsilon_{22}(t) X^\epsilon_2(t) |^2 \right]dt\\
& \quad + \frac{1}{2}\left\langle P^\epsilon(0)x,x\right\rangle  + \frac{1}{2}\int_0^T\left\langle P^\epsilon(t) \sigma^\epsilon, \sigma^\epsilon\right\rangle dt
\end{split}
\end{equation}
where
\begin{align*}
\begin{pmatrix}
\widehat{F}_{11}^\epsilon(t) & \widehat{F}_{12}^\epsilon(t)\\
\widehat{F}_{21}^\epsilon(t) & \widehat{F}_{22}^\epsilon(t)
\end{pmatrix} 
&= \begin{pmatrix}
B_{11}^* P^\epsilon_{11}(t) + B_{21}^* (P^\epsilon_{12}(t))^* & \epsilon B_{11}^* P^\epsilon_{12}(t) + B_{21}^* P^\epsilon_{22}(t)\\
- \left[B_{12}^* P^\epsilon_{11}(t) +  B_{22}^* (P^\epsilon_{12}(t))^*\right] & - \left[ \epsilon B_{12}^* P^\epsilon_{12}(t) + B_{22}^* P^\epsilon_{22}(t)\right]
\end{pmatrix}.
\end{align*}
This implies that
\begin{align*}
J^\epsilon(x_1,x_2; \widehat{u}^\epsilon_1,\widehat{u}^\epsilon_2) = \frac{1}{2}\left\langle P^\epsilon(0)x,x\right\rangle  + \frac{1}{2}\int_0^T\left\langle P^\epsilon(t) \sigma^\epsilon, \sigma^\epsilon\right\rangle dt
\end{align*}
and moreover,
\begin{align*}
&J^\epsilon (x_1,x_2; u_1,u_2) - J^\epsilon(x_1,x_2; \widehat{u}^\epsilon_1,\widehat{u}^\epsilon_2)\\
&= - \frac{1}{2}\E\int_0^T\left[ |u_1(t) - \widehat{F}^\epsilon_{11}(t) X_1(t) - \widehat{F}^\epsilon_{12}(t) X^\epsilon_2(t) |^2 - |u_2(t) - \widehat{F}^\epsilon_{21}(t) X_1(t) - \widehat{F}^\epsilon_{22}(t) X^\epsilon_2(t) |^2 \right]dt.
\end{align*}
From the above, we can see that 
\begin{align*}
&J^\epsilon (x_1,x_2; \widehat{u}^\epsilon_1,u_2) - J^\epsilon(x_1,x_2; \widehat{u}^\epsilon_1,\widehat{u}^\epsilon_2)\\
&= \frac{1}{2}\E\int_0^T|u_2(t) - \widehat{F}^\epsilon_{21}(t) X_1(t) - \widehat{F}^\epsilon_{22}(t) X^\epsilon_2(t) |^2 dt \geq 0.
\end{align*}
This implies that $J^\epsilon (x_1,x_2; \widehat{u}^\epsilon_1,u_2) \geq J^\epsilon(x_1,x_2; \widehat{u}^\epsilon_1,\widehat{u}^\epsilon_2)$ for all $u_2 \in \mathbb{A}_T^{k_2}$. Similarly, we can show that $J^\epsilon (x_1,x_2; u_1,\widehat{u}^\epsilon_2) \leq J^\epsilon(x_1,x_2; \widehat{u}^\epsilon_1,\widehat{u}^\epsilon_2)$ for all $u_1 \in \mathbb{A}_T^{k_1}$. Hence $(\widehat{u}^\epsilon_1,\widehat{u}^\epsilon_2)$ is a feedback saddle point. This also implies that the value is given by $V^\epsilon(x_1,x_2) =  J^\epsilon(x_1,x_2; \widehat{u}^\epsilon_1,\widehat{u}^\epsilon_2)$.

To show uniqueness, suppose that there exists another feedback saddle point and optimal state triple $(\widetilde{u}^\epsilon_1, \widetilde{u}^\epsilon_2, \widetilde{X}^\epsilon)$ where
\begin{align*}
\begin{cases}
\widetilde{u}^\epsilon_1(t) = \widetilde{\phi}_1(t, \widetilde{X}^\epsilon(t),\epsilon), \quad t\in [0,T],\\
\widetilde{u}^\epsilon_2(t) = \widetilde{\phi}_2(t, \widetilde{X}^\epsilon(t), \epsilon), \quad t\in [0,T],
\end{cases}
\end{align*}
where $\widetilde{X}^\epsilon(t) = (\widetilde{X}^\epsilon_1(t),\widetilde{X}^\epsilon_2(t)) =  (X_1(t; x_1, \widetilde{u}^\epsilon_1,\widetilde{u}^\epsilon_2),X^\epsilon_2(t; x_2, \widetilde{u}^\epsilon_1,\widetilde{u}^\epsilon_2))$. To see that both feedback saddle points admit the same value function, we can use the definition of a feedback saddle point to see that
\begin{equation}\label{temmp}
J^\epsilon(x_1,x_2; \widetilde{u}^\epsilon_1,\widetilde{u}^\epsilon_2) \leq J^\epsilon(x_1,x_2; \widetilde{u}^\epsilon_1,\widehat{u}^\epsilon_2)
\leq J^\epsilon(x_1,x_2; \widehat{u}^\epsilon_1,\widehat{u}^\epsilon_2)\leq J^\epsilon(x_1,x_2; \widehat{u}^\epsilon_1,\widetilde{u}^\epsilon_2)\leq J^\epsilon(x_1,x_2; \widetilde{u}^\epsilon_1,\widetilde{u}^\epsilon_2).
\end{equation}
Hence
\begin{equation}\label{temmpp}
J^\epsilon(x_1,x_2; \widetilde{u}^\epsilon_1,\widetilde{u}^\epsilon_2) = J^\epsilon(x_1,x_2; \widehat{u}^\epsilon_1,\widehat{u}^\epsilon_2) = V^\epsilon(x_1,x_2) = \frac{1}{2}\left\langle P^\epsilon(0)x,x\right\rangle  + \frac{1}{2}\int_0^T\left\langle P^\epsilon(t) \sigma^\epsilon, \sigma^\epsilon\right\rangle dt.
\end{equation}
Let us denote $m(t,x',\epsilon) = \widetilde{\phi}_2(t,x',\epsilon) - \widehat{\phi}_2(t,x',\epsilon)$ for $(t,x',\epsilon)\in [0,T]\times \R^n \times (0,1]$. From \eqref{temp0}, \eqref{temmp} and \eqref{temmpp}, we have that
\begin{equation}
\E \int_0^T |m(t, X^{temp,\epsilon}(t),\epsilon)|^2 dt = 0
\end{equation}
where $X^{temp,\epsilon}(t) := X(t;x, \widehat{u}^\epsilon_1, \widetilde{u}^\epsilon_2)$. This implies that
\begin{equation}
m(t, X^{temp,\epsilon}(t),\epsilon) = \widetilde{\phi}_2(t,X^{temp,\epsilon}(t),\epsilon) - \widehat{\phi}_2(t,X^{temp,\epsilon}(t),\epsilon) = 0, \quad \lambda \oplus \PP - a.s.,
\end{equation}
where $\lambda$ is the Lebesgue measure. Moreover, this shows that $X^{temp,\epsilon}$ is a solution to the SDE \eqref{optimstate}. However due to uniqueness of this SDE, we have that $X^{temp,\epsilon} = \widehat{X}^\epsilon$. Thus 
\begin{equation}
\widetilde{\phi}_2(t,\widehat{X}^\epsilon(t),\epsilon) = \widehat{\phi}_2(t,\widehat{X}^\epsilon(t),\epsilon) , \quad \lambda \oplus \PP - a.s.
\end{equation}
for every $\epsilon \in (0,1]$. In the same way, by considering the difference $\widetilde{\phi}_1(t,x',\epsilon) - \widehat{\phi}_1(t,x',\epsilon)$ for $(t,x',\epsilon)\in [0,T]\times \R^n \times (0,1]$, we can show that 
\begin{equation}
\widetilde{\phi}_1(t,\widehat{X}^\epsilon(t),\epsilon) = \widehat{\phi}_1(t,\widehat{X}^\epsilon(t),\epsilon) , \quad \lambda \oplus \PP - a.s.
\end{equation}
for every $\epsilon \in (0,1]$. This implies that $\widetilde{X}^\epsilon$ is also a solution to \eqref{optimstate}. By the uniqueness of the solution to \eqref{optimstate}, we have that $\widetilde{X}^\epsilon= \widehat{X}^\epsilon$. Hence, the triples $(\widehat{u}_1^\epsilon, \widehat{u}_2^\epsilon, \widehat{X}^\epsilon)$ and $(\widetilde{u}_1^\epsilon, \widetilde{u}_2^\epsilon, \widetilde{X}^\epsilon)$ are identical and we arrive at a contradiction.
\qed

\section{Singular perturbation of the Riccati equation}\label{section: reduced}

The purpose of this section is to analyse a version of the generalised Riccati equation when $\epsilon$ is formally set to be $0$, called the reduced system. The main result of this section is to show that reduced system is equivalent to a decoupled pair of differential and algebraic Riccati equations of a lower dimensionality.

Fix $0 < \epsilon \leq 1$. By the following first-order form
\begin{equation}\label{solutionstructure}
P^\epsilon(t) = \begin{pmatrix}
P_{11}^\epsilon(t) & \epsilon P_{12}^\epsilon(t)\\
\epsilon (P_{12}^\epsilon(t))^* & \epsilon P_{22}^\epsilon(t)
\end{pmatrix}, \quad t\in [0,T],
\end{equation}
we can reformulate the generalised Riccati equation \eqref{riccati} in terms of a so-called \textit{full system} of ODEs
\begin{subnumcases}{\label{full}}
\label{pre1}
\frac{dP_{11}^\epsilon}{dt} + f(P_{11}^\epsilon,P_{12}^\epsilon,P_{22}^\epsilon, \epsilon) = 0, \quad P_{11}^\epsilon(T) = 0,\\
\label{pre2}
\epsilon \frac{dP_{12}^\epsilon}{dt} + g_1(P_{11}^\epsilon,P_{12}^\epsilon,P_{22}^\epsilon, \epsilon) = 0, \quad P_{12}^\epsilon(T) = 0,\\
\label{pre3}
\epsilon \frac{dP_{22}^\epsilon}{dt} + g_2(P_{11}^\epsilon,P_{12}^\epsilon,P_{22}^\epsilon, \epsilon) = 0, \quad P_{22}^\epsilon(T) = 0,
\end{subnumcases}
where the functions $f,g_1$ and $g_2$ are defined as
\begin{equation}
\begin{split}
f(P_{11},P_{12},P_{22}, \epsilon) &= A_{11}^* P_{11} + A_{21}^* P_{12}^* + P_{11} A_{11} + P_{12} A_{21} \\
&\quad + P_{11} \Delta_1 P_{11} +P_{12} \Delta^* P_{11} + P_{11} \Delta P_{12}^* + P_{12} \Delta_2 P_{12}^*  + Q_1,\\
\end{split}
\end{equation}
\begin{equation}
\begin{split}
g_1(P_{11},P_{12},P_{22},\epsilon) &= \epsilon A_{11}^* P_{12} + A_{21}^* P_{22} + P_{11} A_{12} + P_{12} A_{22} + \epsilon P_{11} \Delta_1 P_{12} \\
&\quad + \epsilon P_{12} \Delta^* P_{12} + P_{11} \Delta P_{22} + P_{12} \Delta_2 P_{22},\\
\end{split}
\end{equation}
\begin{equation}
\begin{split}
g_2(P_{11},P_{12},P_{22},\epsilon) &= \epsilon A_{12}^* P_{12} + A_{22}^* P_{22} + \epsilon P_{12}^* A_{12} +  P_{22} A_{22} + \epsilon^2 P_{12}^* \Delta_1 P_{12} \\
&\quad + \epsilon P_{12}^* \Delta P_{22} + \epsilon P_{22} \Delta^* P_{12} + P_{22} \Delta_2 P_{22} + Q_2.
\end{split}
\end{equation}
It is clear that we have the following lemma.
\begin{lemma}
For fixed $0< \epsilon \leq 1$, the full system admits a solution $(P_{11}^\epsilon, P_{12}^\epsilon, P_{22}^\epsilon) \in C([0,T]; \mathbb{S}^{n_1}) \times C([0,T]; \R^{n_1\times n_2}) \times C([0,T];\mathbb{S}^{n_2})$ if and only if the generalised Riccati equation \eqref{riccati} admits a solution $P^\epsilon \in C([0,T];\mathbb{S}^n)$.
\end{lemma}

A well-known approach to obtaining the existence, uniqueness and convergence of the solution to the full system is to apply a version of the Tikhonov theorem, see for example Theorem 9.1 of \cite{khalil1996nonlinear}. In order to do so, we need to consider the system of differential and algebraic equations that arise when $\epsilon$ is formally set to be $0$ in the full system. We call this the \textit{reduced system} and it can be represented as the following
\begin{subnumcases}{\label{reducedcompact}}
\label{reduced1}
\frac{d\overline{P}_{11}}{dt} + f(\overline{P}_{11},\overline{P}_{12},\overline{P}_{22}, 0) =0,\quad \overline{P}_{11}(T) = 0,\\
\label{reduced2}
g_1(\overline{P}_{11},\overline{P}_{12},\overline{P}_{22}, 0) = 0,\\
\label{reduced3}
g_2(\overline{P}_{11},\overline{P}_{12},\overline{P}_{22}, 0) = 0.
\end{subnumcases} 
We make the following technical assumption.
\begin{assumption}\label{assumption: tech}
The matrix operator $\Delta_2 = B_{21} B_{21}^* - B_{22} B_{22}^*$ is invertible.
\end{assumption}
Under Assumption \ref{assumption: tech}, let us introduce a decoupled pair of differential and algebraic Riccati equations, which we shall refer to as \textit{reduced differential-algebraic Riccati equations}
\begin{subnumcases}{}
\frac{d\overline{P}_{11}}{dt} + \widetilde{A}^* \overline{P}_{11} + \overline{P}_{11} \widetilde{A} + \overline{P}_{11} M \overline{P}_{11} + N = 0,\quad \overline{P}_{11}(T) = 0, \label{reducedDRE} \\
A_{22}^* \overline{P}_{22} +  \overline{P}_{22} A_{22} +  \overline{P}_{22}\Delta_2 \overline{P}_{22} + Q_2 = 0, \label{reducedARE}
\end{subnumcases}
where
\begin{align*}
\begin{cases}
\widetilde{A} &= A_{11} - \Delta \Delta_2^{-1} \left[ \Delta_2 - A_{22} \Lambda^{-1} A_{22}^* \right] \Delta_2^{-1} A_{21}  - A_{12} \Lambda^{-1} A_{22}^* \Delta_2^{-1} A_{21},\\
M &= \Delta_1 + A_{12} \Lambda^{-1} A_{12}^*- \Delta \Delta_2^{-1} A_{22} \Lambda^{-1} A_{12}^* - A_{12} \Lambda^{-1} A_{22}^* \Delta_2^{-1} \Delta^*  \\
&\quad -\Delta \Delta_2^{-1} \left[  \Delta_2 - A_{22} \Lambda^{-1} A_{22}^*  \right] \Delta_2^{-1} \Delta^*, \\
N &= Q_1 - A_{21}^* \Delta_2^{-1} \left[ \Delta_2 - A_{22}\Lambda^{-1} A_{22}^* \right] \Delta_2^{-1} A_{21}.
\end{cases}
\end{align*}
Here $\Lambda :=  A_{22}^* \Delta_2^{-1} A_{22} - Q_2$. We shall see in Theorem \ref{theorem: reducedRiccati}, that under Assumption \ref{assumption: tech} and \ref{assumption: main}-(b), $\Lambda$ is indeed invertible. We shall refer to \eqref{reducedDRE} as the \textit{reduced differential Riccati equation} and \eqref{reducedARE} as the \textit{reduced algebraic Riccati equation}. We notice that solutions of \eqref{reducedDRE} and \eqref{reducedARE} take values in $\mathbb{S}^{n_1}$ and $\mathbb{S}^{n_2}$, respectively. This is comparatively less than the solutions of the generalised Riccati equation, which take values in the space $\mathbb{S}^{n_1+ n_2}$.

The main assumption of this paper is the existence and uniqueness of a solution to the reduced differential-algebraic Riccati equations \eqref{reducedDRE}-\eqref{reducedARE}.
\begin{assumption}\label{assumption: main}\
\begin{enumerate}[label=(\alph*)]
\item The reduced differential Riccati equation \eqref{reducedDRE} admits a unique solution $\overline{P}^+_{11} \in C([0,T]; \mathbb{S}^{n_1})$. Let $p_0$ be a positive constant such that $|\overline{P}^+_{11}(t)| \leq p_0$ for all $t\in [0,T]$.
\item The reduced algebraic Riccati equation \eqref{reducedARE} admits a stabilising solution $\overline{P}^+_{22} \in \mathbb{S}^{n_2}$. That is, the eigenvalues of $A_{22} + \Delta_2 \overline{P}^+_{22}$ have strictly negative real parts.
\end{enumerate}
\end{assumption}
Note that we do not assume that the reduced algebraic Riccati equation admits a unique solution. However, we require that $\overline{P}^+_{22}$ is an isolated root. This is shown in Proposition 7.9.2 of \cite{lancaster1995algebraic}. We restate this result below.
\begin{proposition}\label{proposition: isolatedroot}
Suppose that Assumption \ref{assumption: main}-(b) holds. Then $\overline{P}^+_{22}$ is a unique solution of \eqref{reducedARE} in the set of solutions $\overline{P}_{22}$ such that the eigenvalues of $A_{22} + \Delta \overline{P}_{22}$ have non-positive real parts.
\end{proposition}

\begin{remark}
Like the original generalised Riccati equation \eqref{riccati}, the reduced differential Riccati equation \eqref{reducedDRE} is also of the generalised type. Sufficient conditions for the existence of a solution to generalised Riccati equation have been formulated to certain extents in \cite{jiongmin2014differential, mcasey2006generalized, yong1999linear,yu2015optimal}. However, analytically checkable sufficient conditions are rather limited. For example, in \cite{yu2015optimal}, the author propose commutative and definiteness conditions, whilst in \cite{jiongmin2014differential}, the author analyses certain examples including the one dimensional case. We will elaborate on the one dimensional case in Section \ref{section: onedimension}. We should also point out that under certain assumptions, the generalised Riccati equation is of the classical type seen in optimal control theory and admits a unique positive/negative definite solution, see Chapter 6 of \cite{yong1999stochastic}. Likewise, under certain assumptions, the algebraic Riccati equation \eqref{reducedARE} admits a unique stablising positive/negative definite solution as seen in optimal control theory, see \cite{lancaster1995algebraic}.
\end{remark}

The following Theorem demonstrates that the reduced system \eqref{reducedcompact} can be characterised in terms of the decoupled pair of reduced differential-algebraic Riccati equations \eqref{reducedDRE}-\eqref{reducedARE}.

\begin{theorem}\label{theorem: reducedRiccati}
Suppose that Assumptions \ref{assumption: tech} and \ref{assumption: main} hold. Let $(\overline{P}^+_{11}, \overline{P}_{22}^+) \in C([0,T];\mathbb{S}^{n_1}) \times \mathbb{S}^{n_2}$ be the solution to the reduced differential-algebraic Riccati equations \eqref{reducedcompact} defined in Assumption \ref{assumption: main}. Then $(\overline{P}^+_{11}, \overline{P}_{12}, \overline{P}_{22}^+) \in C([0,T];\mathbb{S}^{n_1}) \times C([0,T]; \R^{n_1\times n_2}) \times \mathbb{S}^{n_2}$ is a solution to the reduced system \eqref{reducedcompact} where
\begin{equation}
\overline{P}_{12}(t) = -\left(A_{21}^* \overline{P}^+_{22} + \overline{P}^+_{11}(t) A_{12} + \overline{P}^+_{11}(t) \Delta \overline{P}^+_{22}\right)\left( A_{22} + \Delta_2 \overline{P}^+_{22}\right)^{-1}, \quad t\in [0,T].
\end{equation}
\end{theorem}
\noindent \textit{Proof.} Since the algebraic equations \eqref{reduced3} and \eqref{reducedARE} are equivalent, it is clear from Proposition \ref{proposition: isolatedroot} that $\overline{P}^+_{22}$ is the unique stabilising solution of \eqref{reduced3}. From Assumption \ref{assumption: main}-(b), we have that $A_{22} + \Delta_2 \overline{P}^+_{22}$ is invertible. Thus, we can rewrite \eqref{reduced2} as
\begin{equation}\label{temp01}
\overline{P}_{12} = -\left(A_{21}^* \overline{P}^+_{22} + \overline{P}_{11} A_{12} + \overline{P}_{11} \Delta \overline{P}^+_{22}\right)\left( A_{22} + \Delta_2 \overline{P}^+_{22}\right)^{-1}.
\end{equation}
In addition, by Assumption \ref{assumption: tech}, the matrix $\Delta_2$ is invertible and thus, from \eqref{reducedARE} we have that
\[
\Lambda =  A_{22}^* \Delta_2^{-1} A_{22} - Q_2 = (A_{22} + \Delta_2 \overline{P}_{22})^* \Delta_2^{-1} (A_{22} + \Delta_2 \overline{P}_{22})
\]
is invertible. More precisely,
\begin{equation}\label{temp02}
\Lambda^{-1} = (A_{22} + \Delta_2 \overline{P}_{22})^{-1} \Delta_2 (A_{22} + \Delta_2 \overline{P}_{22})^{-*}.
\end{equation}
Applying \eqref{temp01} and \eqref{temp02}, the differential equation \eqref{reduced1} can be written as
\begin{align*}
\frac{d\overline{P}_{11}}{dt} + (A_{11} + I_1)^* \overline{P}_{11}  + \overline{P}_{11} (A_{11} + I_1) + \overline{P}_{11} (\Delta_1 + I_2) \overline{P}_{11} + (Q_1 + I_3)= 0
\end{align*}
where
\begin{align*}
\begin{cases}
I_1 &= - \Delta \Delta_2^{-1} (A_{22} + \Delta_2 \overline{P}_{22}^+) \Lambda^{-1} \overline{P}_{22}^+ A_{21} - (A_{12} + \Delta \overline{P}_{22}^+) \Lambda^{-1} (A_{22} + \Delta_2 \overline{P}_{22}^+)^* \Delta_2^{-1} A_{21}\\
& \quad + (A_{12} + \Delta \overline{P}_{22}^+ ) \Lambda^{-1} \overline{P}_{22}^+ A_{21}\\
I_2 &= - \Delta \Delta_2^{-1} (A_{22} + \Delta_2 \overline{P}_{22}^+) \Lambda^{-1} (A_{12}^* + \overline{P}_{22}^+ \Delta^*) - (A_{12} + \Delta \overline{P}_{22}^+) \Lambda^{-1} (A_{22} + \Delta_2 \overline{P}_{22}^+)^* \Delta_2^{-1} \Delta^*\\
& \quad + (A_{12} + \Delta \overline{P}_{22}^+ ) \Lambda^{-1} ( A_{12}^* + \overline{P}_{22}^+ \Delta^*)\\
I_3 &= -A_{21}^* \Delta_2^{-1} (A_{22} + \Delta_2 \overline{P}_{22}^+) \Lambda^{-1} \overline{P}_{22}^+ A_{21} - A_{21}^* \overline{P}_{22}^+ \Lambda^{-1} (A_{22} + \Delta_2 \overline{P}_{22}^+)^* \Delta_2^{-1} A_{21} \\
&\quad + A_{21}^* \overline{P}_{22}^+ \Lambda^{-1} \overline{P}_{22}^+ A_{21}
\end{cases}
\end{align*}
The rest of the proof involves simplifying the above expressions and removing the dependence on $\overline{P}_{22}^+$. Using \eqref{temp02}, we see that
\begin{align*}
I_1 &= - \Delta \Delta_2^{-1} (A_{22} + \Delta_2 \overline{P}_{22}^+) \Lambda^{-1} \overline{P}_{22}^+ A_{21} - (A_{12} + \Delta \overline{P}_{22}^+) \Lambda^{-1} (A_{22} + \Delta_2 \overline{P}_{22}^+)^* \Delta_2^{-1} A_{21}\\
& \quad + (A_{12} + \Delta \overline{P}_{22}^+ ) \Lambda^{-1} \overline{P}_{22}^+ A_{21}\\
&= - \Delta \Delta_2^{-1} A_{22} \Lambda^{-1} \overline{P}_{22}^+ A_{21} - \Delta \overline{P}_{22}^+ \Lambda^{-1} (A_{22} + \Delta_2 \overline{P}_{22}^+)^* \Delta_2^{-1} A_{21}  - A_{12} \Lambda^{-1} A_{22}^* \Delta_2^{-1} A_{21} \\
&= - \Delta \Delta_2^{-1} \left[ A_{22} \Lambda^{-1} \overline{P}_{22}^+ \Delta_2 + \Delta_2 \overline{P}_{22}^+ \Lambda^{-1} (A_{22} + \Delta_2 \overline{P}_{22}^+)^* \right] \Delta_2^{-1} A_{21}  - A_{12} \Lambda^{-1} A_{22}^* \Delta_2^{-1} A_{21} \\
&= - \Delta \Delta_2^{-1} \left[ A_{22} \Lambda^{-1} (A_{22} + \Delta_2 \overline{P}_{22}^+)^* + \Delta_2 \overline{P}_{22}^+ \Lambda^{-1} (A_{22} + \Delta_2 \overline{P}_{22}^+)^* - A_{22} \Lambda^{-1} A_{22}^* \right] \Delta_2^{-1} A_{21} \\
&\quad - A_{12} \Lambda^{-1} A_{22}^* \Delta_2^{-1} A_{21} \\
&= - \Delta \Delta_2^{-1} \left[ \Delta_2 - A_{22} \Lambda^{-1} A_{22}^* \right] \Delta_2^{-1} A_{21}  - A_{12} \Lambda^{-1} A_{22}^* \Delta_2^{-1} A_{21}.
\end{align*}
Similarly,
\begin{align*}
I_2 &= (A_{12} + \Delta \overline{P}_{22}^+) \Lambda^{-1} (A_{12}^* + \overline{P}_{22}^+ \Delta^*) - \Delta \Delta_2^{-1} (A_{22} + \Delta_2 \overline{P}_{22}^+) \Lambda^{-1} (A_{12}^* + \overline{P}_{22}^+ \Delta^*) \\
& \quad - (A_{12} + \Delta \overline{P}_{22}^+) \Lambda^{-1} (A_{22} + \Delta_2 \overline{P}_{22}^+)^* \Delta_2^{-1} \Delta^*\\
&= A_{12} \Lambda^{-1} A_{12}^* - \Delta \overline{P}_{22}^+ \Lambda^{-1} \overline{P}_{22}^+ \Delta^* - \Delta \Delta_2^{-1} A_{22} \Lambda^{-1} \overline{P}_{22}^+ \Delta^* - \Delta \overline{P}_{22}^+ \Lambda^{-1} A_{22}^* \Delta_2^{-1} \Delta^*\\
&\quad - \Delta \Delta_2^{-1} A_{22} \Lambda^{-1} A_{12}^* - A_{12} \Lambda^{-1} A_{22}^* \Delta_2^{-1} \Delta^*\\
&= A_{12} \Lambda^{-1} A_{12}^*- \Delta \Delta_2^{-1} A_{22} \Lambda^{-1} A_{12}^* - A_{12} \Lambda^{-1} A_{22}^* \Delta_2^{-1} \Delta^*\\
&\quad -\Delta \Delta_2^{-1} \left[  \Delta_2  \overline{P}_{22}^+ \Lambda^{-1} \overline{P}_{22}^+ \Delta_2 + A_{22} \Lambda^{-1} \overline{P}_{22}^+ \Delta_2 + \Delta_2 \overline{P}_{22}^+ \Lambda^{-1} A_{22}^*  \right] \Delta_2^{-1} \Delta^*\\
&= A_{12} \Lambda^{-1} A_{12}^*- \Delta \Delta_2^{-1} A_{22} \Lambda^{-1} A_{12}^* - A_{12} \Lambda^{-1} A_{22}^* \Delta_2^{-1} \Delta^*\\
&\quad -\Delta \Delta_2^{-1} \left[  (A_{22} + \Delta_2  \overline{P}_{22}^+) \Lambda^{-1}(A_{22} + \Delta_2 \overline{P}_{22}^+ )^* - A_{22} \Lambda^{-1} A_{22}^*  \right] \Delta_2^{-1} \Delta^*\\
&= A_{12} \Lambda^{-1} A_{12}^*- \Delta \Delta_2^{-1} A_{22} \Lambda^{-1} A_{12}^* - A_{12} \Lambda^{-1} A_{22}^* \Delta_2^{-1} \Delta^*  -\Delta \Delta_2^{-1} \left[  \Delta_2 - A_{22} \Lambda^{-1} A_{22}^*  \right] \Delta_2^{-1} \Delta^*.
\end{align*}
Lastly,
\begin{align*}
I_3 &= -A_{21}^* \Delta_2^{-1} (A_{22} + \Delta_2 \overline{P}_{22}^+) \Lambda^{-1} \overline{P}_{22}^+ A_{21} - A_{21}^* \overline{P}_{22}^+ \Lambda^{-1} (A_{22} + \Delta_2 \overline{P}_{22}^+)^* \Delta_2^{-1} A_{21} \\
&\quad  + A_{21}^* \overline{P}_{22}^+ \Lambda^{-1} \overline{P}_{22}^+ A_{21}\\
&= -A_{21}^* \Delta_2^{-1} A_{22} \Lambda^{-1} \overline{P}_{22}^+ A_{21} - A_{21}^* \overline{P}_{22}^+ \Lambda^{-1} A_{22}^* \Delta_2^{-1} A_{21}  - A_{21}^* \overline{P}_{22}^+ \Lambda^{-1} \overline{P}_{22}^+ A_{21}\\
&= -A_{21}^* \Delta_2^{-1} \left[ A_{22} \Lambda^{-1} \overline{P}_{22}^+ \Delta_2 + \Delta_2 \overline{P}_{22}^+ \Lambda^{-1} A_{22}^* + \Delta_2 \overline{P}_{22}^+ \Lambda^{-1} \overline{P}_{22}^+ \Delta_2 \right] \Delta_2^{-1} A_{21}\\
&= - A_{21}^* \Delta_2^{-1} \left[ (A_{22} + \Delta_2 \overline{P}_{22}^+ )\Lambda^{-1} (A_{22} + \Delta_2 \overline{P}_{22}^+)^* - A_{22}\Lambda^{-1} A_{22}^* \right] \Delta_2^{-1} A_{21}\\
&= - A_{21}^* \Delta_2^{-1} \left[ \Delta_2 - A_{22}\Lambda^{-1} A_{22}^* \right] \Delta_2^{-1} A_{21}.
\end{align*}
Observe that the differential equation \eqref{reduced1} is equivalent to the reduced differential Riccati equation \eqref{reducedDRE}. Thus, the unique solution $\overline{P}_{11}^+(\cdot)$ to \eqref{reducedDRE} is also the unique solution to \eqref{reduced1}. Finally, from \eqref{temp01}, we deduce that the unique solution of \eqref{reduced2} is indeed given by
\[
\overline{P}_{12}(t) = -\left(A_{21}^* \overline{P}^+_{22} + \overline{P}^+_{11}(t) A_{12} + \overline{P}^+_{11}(t) \Delta \overline{P}^+_{22}\right)\left( A_{22} + \Delta_2 \overline{P}^+_{22}\right)^{-1}, \quad t\in [0,T].
\]
\qed

\section{The Tikhonov Theorem and well-posedness of the generalised Riccati equation}\label{section: asymptotic}

The focus of this section is to establish existence, uniqueness and convergence results for the solution $P^\epsilon$ to the generalised Riccati equation \eqref{riccati}. From Theorem \ref{theorem: reducedRiccati}, we saw that $( \overline{P}_{11},\overline{P}_{12},\overline{P}_{22})$ is a solution to the reduced system. It may be tempting to conclude that the solution of the reduced system is the limit of the solution of the full system, however, not all conditions of the full system are satisfied. Take for example, the terminal condition $P_{22}^\epsilon(T) = 0$. It is very rare that the equality $\lim_{\epsilon \rightarrow 0} P_{22}^\epsilon(T) = 0 = \overline{P}_{22}(T)$ holds. This discrepancy arises because of the degenerate nature of the differential terms and is a featuring characteristic in singular perturbation problems. As we shall see in Theorem \ref{tikhonov}, the reduced system is enough to characterise the limiting solution of the full system over a subinterval $[0,S], S < T$ but not over the entire interval $[0,T]$. Heuristically, this is because the reduced system \eqref{reducedcompact} only captures the fast component of the limit.

To describe the slow component of the limit, we apply a change of time variable $\tau = (T - t)/\epsilon$ and consider the so-called \textit{boundary-layer system}
\begin{subnumcases}{\label{boundarylayer}}
\label{boundarylayer1}
\frac{d\widehat{P}_{12}}{d\tau} = g_1(0,\widehat{P}_{12} + h(0),\widehat{P}_{22} + \overline{P}_{22},0), \quad \widehat{P}_{12}(0) = -h(0),\\
\label{boundarylayer2}
\frac{d\widehat{P}_{22}}{d\tau} = g_2(0,\widehat{P}_{12} + h(0),\widehat{P}_{22} + \overline{P}_{22},0), \quad \widehat{P}_{22}(0) = - \overline{P}_{22}.
\end{subnumcases}
Here the function $h: \mathbb{S}^{n_1} \rightarrow \R^{n_1\times n_2}$ is defined as 
\begin{equation}\label{boundarylayersolution}
h(P_{11}) := - (A_{21}^* \overline{P}_{22} + P_{11} A_{12} + P_{11} \Delta \overline{P}_{22} ) ( A_{22} + \Delta_2 \overline{P}_{22})^{-1}.
\end{equation}
Recall from Proposition \ref{proposition: isolatedroot} and Theorem \ref{theorem: reducedRiccati}, that $h(0)$ and $\overline{P}_{22}$ are an isolated root of 
\[
g_1(0,\overline{P}_{12},\overline{P}_{22},0) = 0\quad \text{ and }\quad g_2(0,\overline{P}_{12},\overline{P}_{22},0) = 0,
\]
respectively. It is clear from this that $(\widehat{P}_{12},\widehat{P}_{22}) = (0,0)$ is an equilibrium of the boundary-layer system. We would like to show this equilibrium is exponentially stable.

Let us introduce some notation. Denote $S = A_{22} + \Delta_2 \overline{P}_{22}$. We strengthen the stability assumption on $S$.
\begin{assumption}\label{assumption: negativedefinite}
The matrix $\frac{1}{2}(S^* + S)$ is negative definite. That is, there exists a positive constant $\gamma$ such that $\frac{1}{2}(S^* + S) \leq -\gamma I$
\end{assumption}
Indeed, under Assumption \ref{assumption: negativedefinite}, the eigenvalues of the matrix $S$ have negative real parts and the stability requirement in Assumption \ref{assumption: main} is satisfied. Define the region 
\[
\mathcal{R}_\delta = \left\{ P_{22} \in \mathbb{S}^{n_2}: \frac{1}{2}\left[ (S + \Delta_2 P_{22})^* + (S + \Delta_2 P_{22}) \right] \leq \delta I \right\}
\]
where $0 < \delta < \gamma$. It is clear that $0 \in \mathcal{R}_\delta$ and for all $P_{22} \in \mathcal{R}_\delta$,
\begin{equation}
\frac{1}{2}\left[ (S^* + S) + (S + \Delta_2 P_{22})^* + (S + \Delta_2 P_{22}) \right] \leq -(\gamma - \delta) I < 0.
\end{equation}
Denote $B_{\delta,q_2}$ as the largest closed ball of radius $q_2$ contained within $\mathcal{R}_\delta$. To ensure the solution converges, we assume that the initial value $-\overline{P}_{22}$ lies in the interior of this ball.
\begin{assumption}\label{assumption: attraction}
For some $ \delta \in (0, \gamma)$, $-\overline{P}_{22}$ is in the interior of closed ball $B_{\delta,q_2}$.
\end{assumption}
Let $\Phi: \mathbb{S}^{n_1} \rightarrow \R^{n_1 \times n_2}$ be the mapping defined as 
\begin{equation}
\Phi(P_{11}) = A_{21}^* + P_{11} \Delta + h(P_{11}) \Delta_2
\end{equation}
and let $q_1$ be a positive constant such that $q_1 > e^{\frac{|\Delta_2| q_2}{\gamma - \delta}} \left(|h(0)| + \frac{|\Phi(0)| q_2}{\delta}\right) + |h(0)|$. We denote $B_{q_1}$ as the closed ball in $\R^{n_1\times n_2}$ of radius $q_1$.
\begin{lemma}\label{lemma: exponentialstabilising}
Suppose that Assumptions \ref{assumption: tech}, \ref{assumption: main}, \ref{assumption: negativedefinite} and \ref{assumption: attraction} hold. Then there exists a unique solution $(\widehat{P}_{12}(\tau),\widehat{P}_{22}(\tau))$ to the boundary-layer system \eqref{boundarylayer} contained in $B_{q_1} \times B_{\delta,q_2}$, for all $\tau \geq 0$, which converges to $(0,0)$. Moreover, the equilibrium $(0,0)$ is exponentially stable with an estimated region of attraction $R_A = \{(P_{12},P_{22}) : P_{22} \in B_{\delta,q_2}\}$. To be precise, there exists positive constants $k_1,k_2$, which may depend on $q_2,\gamma,\delta$, such that for all initial values $(\widehat{P}_{12}(0),\widehat{P}_{22}(0)) \in R_A$,
\begin{equation}\label{stab1}
|\widehat{P}_{12}(\tau)| \leq  k_1 e^{-\gamma \tau} |\widehat{P}_{12}(0)| + k_2 e^{-(\gamma -\delta) \tau}, \quad  \forall \tau \geq 0,
\end{equation}
and
\begin{equation}\label{stab2}
|\widehat{P}_{22}(\tau)| \leq e^{-(\gamma - \delta) \tau} |\widehat{P}_{22}(0)|, \quad \forall \tau \geq 0.
\end{equation}
\end{lemma}
\noindent \textit{Proof.} Evaluating the boundary-layer system \eqref{boundarylayer}, we have that
\begin{subnumcases}{}
\label{blay1}
\frac{d\widehat{P}_{12} }{d\tau} = \widehat{P}_{12} S + \Phi(0) \widehat{P}_{22} + \widehat{P}_{12} \Delta_2 \widehat{P}_{22},\\
\label{blay2}
\frac{d\widehat{P}_{22}}{d\tau} = S^* \widehat{P}_{22} + \widehat{P}_{22} S + \widehat{P}_{22} \Delta_2 \widehat{P}_{22}.
\end{subnumcases}
Starting with \eqref{blay2}, differentiating the squared norm of $\widehat{P}_{22}$, we have that
\begin{align*}
\frac{d}{d\tau} |\widehat{P}_{22}|^2 &= 2 \left\langle \frac{d\widehat{P}_{22} }{d\tau}, \widehat{P}_{22} \right\rangle\\
&= 2 \langle S^* \widehat{P}_{22} + \widehat{P}_{22} S + \widehat{P}_{22} \Delta_2 \widehat{P}_{22} , \widehat{P}_{22} \rangle\\
&= \langle S^* \widehat{P}_{22} + \widehat{P}_{22} S, \widehat{P}_{22} \rangle + \langle (S + \Delta_2 \widehat{P}_{22})^* \widehat{P}_{22} + \widehat{P}_{22} (S + \Delta_2 \widehat{P}_{22} ) , \widehat{P}_{22} \rangle.
\end{align*}
Thus, for $\widehat{P}_{22} \in B_{\delta, q_2}$,
\begin{equation}\label{temp-1}
\frac{d}{d\tau} |\widehat{P}_{22}(\tau)|^2 \leq  - 2 (\gamma - \delta) |\widehat{P}_{22}(\tau)|^2
\end{equation}
which implies that $\widehat{P}_{22}(\tau)$ is exponentially stable
\begin{equation}\label{temp-2}
|\widehat{P}_{22}(\tau)| \leq e^{-(\gamma - \delta) \tau} |\widehat{P}_{22}(0)|,\quad \forall \tau \geq 0.
\end{equation}
Since $\widehat{P}_{22}(0) = -\overline{P}_{22}$ is contained in the interior of the closed ball $B_{\delta,q_2}$, the inequality \eqref{temp-2} implies that the continuous solutions $\widehat{P}_{22}(\tau)$ of \eqref{boundarylayersolution} are also contained in $B_{\delta,q_2}$. Thus, the solution $\widehat{P}_{22}(\tau)$ exists for all $\tau \geq 0$.  

Similarly for \eqref{blay1}, by the Cauchy-Schwartz inequality, we have that
\begin{align*}
\frac{d}{d\tau}|\widehat{P}_{12}|^2
&= 2 \langle \widehat{P}_{12} S, \widehat{P}_{12} \rangle +  2 \langle \widehat{P}_{12} \Delta_2 \widehat{P}_{22} , \widehat{P}_{12} \rangle + 2 \langle \widehat{P}_{12}, \Phi(0) \widehat{P}_{22} \rangle\\
&\leq - 2\gamma |\widehat{P}_{12}|^2 + 2 |\Delta_2| |\widehat{P}_{22}| |\widehat{P}_{12}|^2  + 2 |\Phi(0)| |\widehat{P}_{22}| |\widehat{P}_{12}|. 
\end{align*}
So for all $\widehat{P}_{22} \in B_{\delta,q_2}$, we can apply \eqref{temp-2} to obtain
\begin{align*}
\frac{d}{d\tau}|\widehat{P}_{12}|^2 &\leq \left( 2 |\Delta_2| e^{-(\gamma - \delta) \tau} |\widehat{P}_{22}(0)| - 2\gamma  \right) |\widehat{P}_{12}|^2  + 2 |\Phi(0)|  e^{-(\gamma - \delta) \tau} |\widehat{P}_{22}(0)| |\widehat{P}_{12}|\\
&\leq \left( 2  q_2|\Delta_2|  e^{-(\gamma - \delta) \tau} - 2 \gamma  \right) |\widehat{P}_{12}|^2  + 2 |\Phi(0)|   q_2 e^{-(\gamma - \delta) \tau}  |\widehat{P}_{12}|.
\end{align*}
Using the chain rule, we can show that $ \frac{d}{d\tau}|\widehat{P}_{12}|^2 = 2|\widehat{P}_{12}| \frac{d}{d\tau}|\widehat{P}_{12}|$. This implies that
\begin{align*}
\frac{d}{d\tau}|\widehat{P}_{12}| \leq \left( q_2|\Delta_2|  e^{-(\gamma - \delta) \tau} - \gamma  \right) |\widehat{P}_{12}|  +  |\Phi(0)|   q_2 e^{-(\gamma - \delta) \tau}.
\end{align*}
By a change of variations approach, we can show the above gives
\begin{align*}
|\widehat{P}_{12}(\tau)| &\leq e^{\frac{q_2 |\Delta_2| }{\gamma - \delta} - \gamma \tau}|\widehat{P}_{12}(0)| + \frac{ |\Phi(0)| q_2}{\delta} e^{\frac{q_2 |\Delta_2| }{\gamma - \delta} - (\gamma - \delta) \tau}\\
&= k_1 e^{-\gamma \tau} |\widehat{P}_{12}(0)| + k_2 e^{-(\gamma -\delta) \tau}
\end{align*}
where $k_1 = e^{\frac{|\Delta_2| q_2 }{\gamma - \delta}}$ and $k_2 = \frac{ |\Phi(0)| q_2}{\delta} e^{\frac{|\Delta_2| q_2 }{\gamma - \delta}}$.
Hence, for all $\tau \geq 0$, the continuous solutions $\widehat{P}_{12}(\tau)$ of \eqref{boundarylayer} are contained in the closed ball $B_{q_1}$ of radius $q_1$. This implies the solution $\widehat{P}_{12}(\tau)$ exists for all $\tau \geq 0$. The uniqueness of the solution to the boundary-layer system \eqref{boundarylayer} follows from the locally Lipschitz property of the functions $g_1$ and $g_2$.
\qed 

The following theorem is the main result of this paper. 

\begin{theorem}\label{tikhonov} 
Suppose that Assumptions \ref{assumption: tech}, \ref{assumption: main}, \ref{assumption: negativedefinite} and \ref{assumption: attraction} hold. Let $T>0$ be any finite time horizon. Then there exists a positive constant $\epsilon^*$ such that, for all $0 < \epsilon < \epsilon^*$, the full system \eqref{full} admits a unique solution $(P^\epsilon_{11}, P^\epsilon_{12}, P^\epsilon_{22}) \in C([0,T];\mathbb{S}^{n_1}) \times C([0,T];\R^{n\times m}) \times \mathbb{S}^{n_2}$, which satisfy the estimates
\begin{equation}\label{tikhonoveq}
\begin{split}
P_{11}^\epsilon(t) - \overline{P}_{11}(t) &= O(\epsilon),\\
P_{12}^\epsilon(t) - \overline{P}_{12}(t) - \widehat{P}_{12}\left(\frac{T-t}{\epsilon}\right) &= O(\epsilon),\\
P_{22}^\epsilon(t) - \overline{P}_{22} - \widehat{P}_{22}\left(\frac{T-t}{\epsilon}\right) &= O(\epsilon),
\end{split}
\end{equation}
uniformly in $t\in [0,T]$. Here $(\overline{P}_{11},\overline{P}_{12},\overline{P}_{22})$ is the solution to the reduced system \eqref{reducedcompact} defined in Theorem \ref{theorem: reducedRiccati} and $(\widehat{P}_{12}(\tau),\widehat{P}_{22}(\tau))$ is the solution of the boundary-layer system \eqref{boundarylayer}.
\end{theorem}
The proof proceeds in a similar manner to the proof of the Tikhonov theorem in \cite{khalil1996nonlinear}. The main difference is that in \cite{khalil1996nonlinear}, the authors apply results from Lyapunov theory and in doing so, require a stronger assumptions on the terminal values of the full system. More specifically, it is possible $0$ is in a smaller neighbourhood of $-\overline{P}_{22}$ than was specified in Assumption \ref{assumption: attraction}. In this paper, we use a direct approach to show that Assumption \ref{assumption: attraction} is indeed sufficient. Section \ref{section: proof} of this paper is dedicated to the proof.

The following corollary describes the existence, uniqueness and convergence of the solution of the generalised Riccati equation \eqref{riccati}.
\begin{theorem}\label{theorem: maincorollary}
Suppose that Assumptions \ref{assumption: tech}, \ref{assumption: main}, \ref{assumption: negativedefinite} and \ref{assumption: attraction} hold. Let $T > 0$ be any finite time horizon and $\epsilon^*$ be the small positive parameter defined in Theorem \ref{tikhonov}. If $0 < \epsilon < \epsilon^*$ then the Riccati equation \eqref{riccati} admits a unique solution
\[
P^\epsilon(t) = \begin{pmatrix}
P_{11}^\epsilon(t) & \epsilon P_{12}^\epsilon(t)\\
\epsilon (P_{12}^\epsilon(t))^* & \epsilon P_{22}^\epsilon(t)
\end{pmatrix}, \quad t\in [0,T],
\]
in the space $C([0,T];\mathbb{S}^n)$, where $(P_{11}^\epsilon(t), P_{12}^\epsilon(t), P_{22}^\epsilon(t))$ are the unique solutions of the full system \eqref{full} described in Theorem \ref{tikhonov}. Moreover, $P^\epsilon$ possesses the asymptotic property
\[
 P^\epsilon(t)- \begin{pmatrix}
\overline{P}_{11}(t) & 0\\
0 & 0
\end{pmatrix} = O(\epsilon), \quad \forall t\in [0,T],
\]
where $\overline{P}_{11}(\cdot)$ is the unique solution of the reduced differential Riccati equation \eqref{reducedDRE}. 
\end{theorem}
\noindent \textit{Proof.} For $0 < \epsilon < \epsilon^*$, the existence of a solution to generalised Riccati equation \eqref{riccati} follows from Theorem \ref{tikhonov} and the first-order form \eqref{solutionstructure}. The uniqueness follows from the locally Lipschitz property of the non-linear terms. We are left to prove the convergence result. For all $i = 1,2$, $t \in [0,T]$ and $0 < \epsilon < \epsilon^*$, Lemma \ref{lemma: exponentialstabilising} implies that $\epsilon  \widehat{P}_{i2}\left( \frac{T - t}{\epsilon} \right)  = O(\epsilon)$. Hence by Theorem \ref{tikhonov}, for all $i = 1,2$, $t \in [0,T]$ and $0 < \epsilon < \epsilon^*$, we have
\[
\epsilon P_{i2}^\epsilon(t) = O(\epsilon).
\]
This completes the proof.
\qed

The following corollary to Theorem \ref{tikhonov} gives a useful estimate for the next section. 

\begin{corollary}\label{corollary: discbounds}
Suppose that Assumptions \ref{assumption: tech}, \ref{assumption: main}, \ref{assumption: negativedefinite} and \ref{assumption: attraction} hold. Let $T > 0$ be any finite time horizon and $\epsilon^*$ be the small positive parameter defined in Theorem \ref{tikhonov}. Then for any positive integer $j$, there exists a positive constant $K(T,j)$, which depends on $T$ and $j$, such that for all $0 < \epsilon < \epsilon^*$ and $i= 1,2$
\begin{equation}\label{discbounds1}
\int_0^T |P_{i2}^\epsilon(t) - \overline{P}_{i2}(t)|^j dt \leq \epsilon K(T,j).
\end{equation}
\end{corollary}
\noindent \textit{Proof.} Fix $T > 0$ and let $t\in [0,T]$. For all $i = 1,2$, Lemma \ref{lemma: exponentialstabilising} and Theorem \ref{tikhonov} imply that
\begin{align*}
| P_{i2}^\epsilon(t) - \overline{P}_{i2}(t) | &\leq \epsilon K_1(T) + \Big| \widehat{P}_{i2}\left(\frac{T - t}{\epsilon}\right) \Big|\\
& \leq \epsilon K_1(T) + K_2 e^{-\frac{(\gamma - \delta)(T - t)}{\epsilon}}
\end{align*}
for some positive constants $K_1(T)$, which depends on $T$, and $K_2$. Hence, for all positive integers $j$
\begin{align*}
| P_{i2}^\epsilon(t) - \overline{P}_{i2}(t) |^j
& \leq 2^{j-1} \left(\epsilon^j K_1(T)^j + K_2^j e^{-\frac{j (\gamma - \delta)(T-t)}{\epsilon}}\right).
\end{align*}
The result follows from integrating the above expression. 
\qed

\section{Approximating saddle point and estimation of the value function}\label{section: approxcontrol}
In this section, we propose an approximate feedback saddle point to Problem (SLQG) based on the reduced system and demonstrate that the objective function associated with this pair of strategies convergences to the value function with a rate of order $O(\epsilon)$. We preface that we denote $K$ as a positive constant and are not necessarily the same in each instance. In situations where $K$ may depend on another relevant constant, say $T$, then we will denote this as $K(T)$.

From Theorem \ref{theorem: solution}, the feedback saddle point of Problem (SLQG) is given by
\begin{equation}\label{splitoptimcontrol}
\begin{split}
\begin{cases}
\widehat{u}^\epsilon_1(t) = \widehat{\phi}_1(t,\widehat{X}^\epsilon_1(t),\widehat{X}^\epsilon_2(t),\epsilon) = \widehat{F}_{11}^\epsilon(t) \widehat{X}^\epsilon_1(t) + \widehat{F}_{12}^\epsilon(t) \widehat{X}^\epsilon_2 (t),\quad t\in [0,T],\\
\widehat{u}^\epsilon_2(t) = \widehat{\phi}_2(t,\widehat{X}^\epsilon_1(t),\widehat{X}^\epsilon_2(t),\epsilon) =  \widehat{F}_{21}^\epsilon(t) \widehat{X}^\epsilon_1(t) + \widehat{F}_{22}^\epsilon(t) \widehat{X}^\epsilon_2 (t),\quad t\in [0,T],
\end{cases}
\end{split}
\end{equation}
where
\begin{align*}
\begin{pmatrix}
\widehat{F}_{11}^\epsilon(t) & \widehat{F}_{12}^\epsilon(t)\\
\widehat{F}_{21}^\epsilon(t) & \widehat{F}_{22}^\epsilon(t)
\end{pmatrix} 
&= \begin{pmatrix}
B_{11}^* P^\epsilon_{11}(t) + B_{21}^* (P^\epsilon_{12}(t))^* & \epsilon B_{11}^* P^\epsilon_{12}(t) + B_{21}^* P^\epsilon_{22}(t)\\
- \left[B_{12}^* P^\epsilon_{11}(t) +  B_{22}^* (P^\epsilon_{12}(t))^*\right] & - \left[ \epsilon B_{12}^* P^\epsilon_{12}(t) + B_{22}^* P^\epsilon_{22}(t)\right]
\end{pmatrix}
\end{align*}
are feedback operators and $\widehat{X}^\epsilon_1$ and $\widehat{X}^\epsilon_2$ are the solution of the state equations
\begin{equation}\label{slowoptim}
\begin{split}
\begin{cases}
d\widehat{X}^\epsilon_1(t) = \left[\left(A_{11} + \Delta_1 P_{11}^\epsilon + \Delta (P_{12}^\epsilon)^*\right) \widehat{X}^\epsilon_1(t) + \left( A_{12}+ \epsilon \Delta_1 P_{12}^\epsilon + \Delta P_{22}^\epsilon \right) \widehat{X}^\epsilon_2 (t)\right]dt \\
\qquad \qquad + \sigma_1 dW_1(t),\\
\widehat{X}^\epsilon_1(0) = x_1,
\end{cases}
\end{split}
\end{equation}
and
\begin{equation}\label{fastoptim}
\begin{split}
\begin{cases}
d\widehat{X}^\epsilon_2 (t) = \frac{1}{\epsilon} \left[\left(A_{21} + \Delta^* P_{11}^\epsilon + \Delta_2 (P_{12}^\epsilon)^* \right)\widehat{X}^\epsilon_1(t) + \left(A_{22} + \epsilon \Delta^* P_{12}^\epsilon + \Delta_2 P_{22}^\epsilon \right) \widehat{X}^\epsilon_2 (t) \right] dt\\
\qquad \qquad + \frac{1}{\sqrt{\epsilon}} \sigma_2 dW_2(t),\\
\widehat{X}^\epsilon_2(0) = x_2.
\end{cases}
\end{split}
\end{equation}
We construct an approximate feedback saddle point $\overline{u}^\epsilon(\cdot) = (\overline{u}^\epsilon_1(\cdot),\overline{u}^\epsilon_2(\cdot))$ by formally setting $\epsilon = 0$ in the feedback operator of \eqref{splitoptimcontrol}
\begin{equation}\label{approxoptimcontrol}
\begin{split}
\begin{cases}
\overline{u}^\epsilon_1(t) := \widehat{\phi}_1(t,\overline{X}^\epsilon_1(t),\overline{X}^\epsilon_2(t),0) = \overline{F}_{11}(t) \overline{X}^\epsilon_1(t) + \overline{F}_{12} \overline{X}^\epsilon_2 (t),\quad t\in [0,T],\\
\overline{u}^\epsilon_2(t) := \widehat{\phi}_2(t,\overline{X}^\epsilon_1(t),\overline{X}^\epsilon_2(t),0) = \overline{F}_{21}(t) \overline{X}^\epsilon_1(t) + \overline{F}_{22} \overline{X}^\epsilon_2 (t),\quad t\in [0,T],
\end{cases}
\end{split}
\end{equation}
where
\begin{align*}
\begin{pmatrix}
\overline{F}_{11}(t) & \overline{F}_{12}\\
\overline{F}_{21}(t) & \overline{F}_{22}
\end{pmatrix} 
= \begin{pmatrix}
B_{11}^* \overline{P}_{11}(t) + B_{21}^* (\overline{P}_{12}(t))^* &  B_{21}^* \overline{P}_{22} \\
- \left[B_{12}^* \overline{P}_{11}(t) +  B_{22}^* (\overline{P}_{12}(t))^*\right] & -  B_{22}^* \overline{P}_{22}
\end{pmatrix}
\end{align*}
and $\overline{X}^\epsilon_1$ and $\overline{X}^\epsilon_2$ are the solution of the following state equations
\begin{equation}\label{slowoptimapprox}
\begin{split}
\begin{cases}
d\overline{X}^\epsilon_1(t) = \left[\left(A_{11} + \Delta_1 \overline{P}_{11}(t) + \Delta \overline{P}_{12}(t)^* \right) \overline{X}^\epsilon_1(t) +  (A_{12} + \Delta \overline{P}_{22}) \overline{X}^\epsilon_2(t)\right]dt + \sigma_1 dW_1(t),\\
\overline{X}^\epsilon_1(0) = x_1,
\end{cases}
\end{split}
\end{equation}
and
\begin{equation}\label{fastoptimapprox}
\begin{split}
\begin{cases}
d \overline{X}^\epsilon_2(t) = \frac{1}{\epsilon} \left[\left(A_{21} + \Delta^* \overline{P}_{11}(t) + \Delta_2 \overline{P}_{12}(t)^* \right) \overline{X}^\epsilon_1(t) + \left(A_{22} + \Delta_2 \overline{P}_{22} \right) \overline{X}^\epsilon_2 (t) \right] dt + \frac{1}{\sqrt{\epsilon}} \sigma_2 dW_2(t) , \\ \overline{X}^\epsilon_2(0) = x_2.
\end{cases}
\end{split}
\end{equation} 
In the following Lemma, we show that the feedback operators of $(\widehat{u}^\epsilon_1, \widehat{u}^\epsilon_2)$ and $(\overline{u}_1^\epsilon,\overline{u}_2^\epsilon)$ converge in $L^2$.

\begin{lemma}\label{lemma: feedbackopertorconverge}
Suppose that Assumptions \ref{assumption: tech}, \ref{assumption: main}, \ref{assumption: negativedefinite} and \ref{assumption: attraction} hold. Let $T > 0$ be any finite time horizon and $\epsilon^*$ be the small positive parameter defined in Theorem \ref{tikhonov}. Then there exists a positive constant $K(T)$, which depends on $T$, such that for all $i,j = 1,2$ and $0 < \epsilon < \epsilon^*$
\begin{equation}
\int_0^T |\widehat{F}_{ij}^\epsilon(t) - \overline{F}_{ij}(t)|^2 dt \leq \epsilon K(T).
\end{equation}
\end{lemma}
\noindent \textit{Proof.} The result follows from Theorem \ref{theorem: maincorollary} and Corollary \ref{corollary: discbounds}.
\qed

\begin{lemma}\label{lemma: boundedtraj}
Suppose that Assumptions \ref{assumption: tech}, \ref{assumption: main}, \ref{assumption: negativedefinite} and \ref{assumption: attraction} hold. Then for any finite $T>0$
\begin{equation}
\sup_{\epsilon\in (0,1]} \sup_{t\in [0,T]} \E\left[ |\overline{X}^\epsilon_1(t) |^2 + |\overline{X}^\epsilon_2(t) |^2\right] < \infty.
\end{equation}
\end{lemma}
\noindent \textit{Proof.} Fix $\epsilon \in (0,1]$. From the slow state process \eqref{slowoptimapprox}, for all $t\in [0,T]$,
\begin{align*}
|\overline{X}^\epsilon_1(t) | &\leq |x_1|  +\int_0^t \left[|A_{11} + \Delta_1 \overline{P}_{11}(s) + \Delta (\overline{P}_{12}(s))^* | \ |\overline{X}^\epsilon_1(s)| + | A_{12} + \Delta \overline{P}_{22} | \ | \overline{X}^\epsilon_2 (s) | \right]ds \\
&\quad + |\sigma_1 | |W_1(t)|. 
\end{align*}
By Theorem \ref{theorem: reducedRiccati}, the processes $\overline{P}_{11}(t)$ and $\overline{P}_{12}(t)$ are uniformly bounded over the interval $t\in [0,T]$. Thus, we have that
\begin{align*}
|\overline{X}^\epsilon_1(t) |\leq |x_1| + K(T) \int_0^t \left[ |\overline{X}^\epsilon_1(s)| + |\overline{X}^\epsilon_2 (s) | \right]ds + |\sigma_1 | |W_1(t)|.
\end{align*}
Squaring the above inequality, and applying Ito's Isometry, we can bound the following expectation
\begin{equation}\label{xbound}
\E \left[|\overline{X}^\epsilon_1(t) |^2 \right] \leq 3|x_1|^2 + 3|\sigma_1|^2 nt  + K(T) \int_0^t\E \left[ | \overline{X}^\epsilon_1(s)|^2 + |\overline{X}^\epsilon_2 (s) |^2\right]ds.
\end{equation}

Now we turn to the fast state process \eqref{fastoptimapprox}. Fix $t \in [0,T]$. Applying Ito's lemma to the mapping $s \mapsto e^{(A_{22} + \Delta_2 \overline{P}_{22} ) \frac{t-s}{\epsilon}} \overline{X}^\epsilon_2(s)$, we have that
\begin{align*}
&d\left( e^{(A_{22} + \Delta_2 \overline{P}_{22} ) \frac{t-s}{\epsilon}} \overline{X}^\epsilon_2(s) \right) \\
&=- \frac{1}{\epsilon}(A_{22} + \Delta_2 \overline{P}_{22} ) e^{(A_{22} + \Delta_2 \overline{P}_{22} ) \frac{t-s}{\epsilon}} \overline{X}^\epsilon_2(s) ds + e^{(A_{22} + \Delta_2 \overline{P}_{22} ) \frac{t-s}{\epsilon}} d\overline{X}^\epsilon_2(s)\\
&= \frac{1}{\epsilon} e^{(A_{22} + \Delta_2 \overline{P}_{22} ) \frac{t-s}{\epsilon}} \left[\left(A_{21} + \Delta^* \overline{P}_{11}(s) + \Delta_2 (\overline{P}_{12}(s))^* \right)\overline{X}^\epsilon_1(s)\right] ds + \frac{1}{\sqrt{\epsilon}} e^{(A_{22} + \Delta_2 \overline{P}_{22} ) \frac{t-s}{\epsilon}}  \sigma_2 dW_2(s).
\end{align*}
Thus, by another application of Ito's lemma to the mapping $s\mapsto |e^{(A_{22} + \Delta_2 \overline{P}_{22} ) \frac{t-s}{\epsilon}}  \overline{X}^\epsilon_2(s)|^2$
\begin{align*}
d & |e^{(A_{22} + \Delta_2 \overline{P}_{22} ) \frac{t-s}{\epsilon}}  \overline{X}^\epsilon_2(s)|^2 \\
&= d \left\langle e^{(A_{22} + \Delta_2 \overline{P}_{22} ) \frac{t-s}{\epsilon}} \overline{X}^\epsilon_2(s) , e^{(A_{22} + \Delta_2 \overline{P}_{22} ) \frac{t-s}{\epsilon}} \overline{X}^\epsilon_2(s) \right\rangle\\
&= 2  \left\langle e^{(A_{22} + \Delta_2 \overline{P}_{22} ) \frac{t-s}{\epsilon}} \overline{X}^\epsilon_2(s) , d \left( e^{(A_{22} + \Delta_2 \overline{P}_{22} ) \frac{t-s}{\epsilon}} \overline{X}^\epsilon_2(s)\right) \right\rangle \\
&\quad + \left\langle d \left( e^{(A_{22} + \Delta_2 \overline{P}_{22} ) \frac{t-s}{\epsilon}} \overline{X}^\epsilon_2(s) \right), d \left( e^{(A_{22} + \Delta_2 \overline{P}_{22} ) \frac{t-s}{\epsilon}} \overline{X}^\epsilon_2(s)\right) \right\rangle\\
&= \frac{1}{\epsilon} |e^{(A_{22} + \Delta_2 \overline{P}_{22} ) \frac{t-s}{\epsilon}}  \sigma_2 |^2 ds + \frac{2}{\sqrt{\epsilon}} \left\langle e^{(A_{22} + \Delta_2 \overline{P}_{22} ) \frac{t-s}{\epsilon}} \overline{X}^\epsilon_2(s), e^{(A_{22} + \Delta_2 \overline{P}_{22} ) \frac{t-s}{\epsilon}}  \sigma_2 dW_2(s) \right\rangle\\
&\quad + \frac{2}{\epsilon} \left\langle e^{(A_{22} + \Delta_2 \overline{P}_{22} ) \frac{t-s}{\epsilon}} \overline{X}^\epsilon_2(s), e^{(A_{22} + \Delta_2 \overline{P}_{22} ) \frac{t-s}{\epsilon}}  \left(A_{21} + \Delta^* \overline{P}_{11}(s) + \Delta_2 (\overline{P}_{12}(s))^* \right)\overline{X}^\epsilon_1(s)\right\rangle ds.
\end{align*}
Integrating from $0$ to $t$ and taking the expectation, we have that
\begin{align*}
&\E \left[ |\overline{X}^\epsilon_2(t)|^2 \right] = |e^{(A_{22} + \Delta_2 \overline{P}_{22} ) \frac{t}{\epsilon}} x_2|^2 + \frac{1}{\epsilon} \int_0^t |e^{(A_{22} + \Delta_2 \overline{P}_{22} ) \frac{t-s}{\epsilon}}  \sigma_2 |^2 ds\\
&\quad + \frac{2}{\epsilon} \E \int_0^t \left\langle e^{(A_{22} + \Delta_2 \overline{P}_{22} ) \frac{t-s}{\epsilon}} \overline{X}^\epsilon_2(s), e^{(A_{22} + \Delta_2 \overline{P}_{22} ) \frac{t-s}{\epsilon}}  \left(A_{21} + \Delta^* \overline{P}_{11}(s) + \Delta_2 (\overline{P}_{12}(s))^* \right)\overline{X}^\epsilon_1(s)\right\rangle ds\\
&\quad + \frac{2}{\sqrt{\epsilon}}\E\int_0^t \left\langle e^{(A_{22} + \Delta_2 \overline{P}_{22} ) \frac{t-s}{\epsilon}} \overline{X}^\epsilon_2(s), e^{(A_{22} + \Delta_2 \overline{P}_{22} ) \frac{t-s}{\epsilon}}  \sigma_2 dW_2(s) \right\rangle.
\end{align*}
Recall from Theorem \ref{theorem: reducedRiccati} that $\overline{P}_{11}(t)$ and $\overline{P}_{12}(t)$ are uniformly bounded on the interval $t\in [0,T]$. Thus, for fixed $\epsilon \in (0,1]$, $\overline{X}^\epsilon_2$ is the solution of a linear stochastic differential equation with drift and diffusion terms satisfying Lipscthiz and linear growth conditions. As result, $\overline{X}^\epsilon_2$ is square integrable for fixed $\epsilon\in (0,1]$ (see for example Theorem 5.2.1 of \cite{oksendal2003stochastic}) and subsequently,
\begin{equation}\label{stochasticintegral}
\frac{2}{\sqrt{\epsilon}} \E \int_0^t \left\langle e^{(A_{22} + \Delta_2 \overline{P}_{22} ) \frac{t-s}{\epsilon}} \overline{X}^\epsilon_2(s), e^{(A_{22} + \Delta_2 \overline{P}_{22} ) \frac{t-s}{\epsilon}}  \sigma_2 dW_2(s) \right\rangle = 0.
\end{equation}
From Assumption \ref{assumption: main}, $A_{22} + \Delta_2 \overline{P}_{22}$ has eigenvalues with negative real parts. Thus there exists positive constants $M$ and $\gamma$ such that
\begin{equation}\label{inequalityexp}
|e^{\left( A_{22} + \Delta_2 \overline{P}_{22} \right) \frac{t}{\epsilon}}| \leq M e^{-\frac{\gamma t}{\epsilon}}, \qquad \forall t\in [0,T].
\end{equation}
Using \eqref{stochasticintegral}, \eqref{inequalityexp}, the uniform boundedness of $\overline{P}_{11}$ and $\overline{P}_{12}$, and the Cauchy-Schwartz inequality, we have that
\begin{align*}
&\E \left[ |\overline{X}^\epsilon_2(t)|^2 \right]  \\
&\leq M^2 e^{-\frac{2\gamma t}{\epsilon}} |x_2|^2 + \frac{M^2 | \sigma_2|^2}{\epsilon} \int_0^t e^{-\frac{2\gamma (t - s)}{\epsilon}} ds + \frac{K(T)}{\epsilon} \int_0^t  e^{-\frac{2\gamma (t - s)}{\epsilon}} \E \left[ |\overline{X}^\epsilon_1(s)|^2 + |\overline{X}^\epsilon_2(s)|^2 \right] ds\\
&= M^2 e^{-\frac{2\gamma t}{\epsilon}} |x_2|^2 + \frac{M^2 |\sigma_2|^2}{2\gamma} (1 - e^{-\frac{2\gamma t}{\epsilon}}) + \frac{K(T)}{\epsilon} \int_0^t  e^{-\frac{2\gamma (t - s)}{\epsilon}} \E \left[ |\overline{X}^\epsilon_1(s)|^2 + |\overline{X}^\epsilon_2(s)|^2 \right] ds\\
&\leq  M^2 |x_2|^2 + \frac{M^2 |\sigma_2|^2}{2\gamma} + \frac{K(T)}{\epsilon} \int_0^t  e^{-\frac{2\gamma (t - s)}{\epsilon}} \E \left[ |\overline{X}^\epsilon_1(s)|^2 + |\overline{X}^\epsilon_2(s)|^2 \right] ds
\end{align*}
Summing the above expression with \eqref{xbound}, we have that
\begin{equation}\label{postgronwall}
\begin{split}
&\E \left[ |\overline{X}^\epsilon_1(t)|^2 +  |\overline{X}^\epsilon_2(t)|^2 \right]\\
&\leq 4|x_1|^2 + 4|\sigma_1|^2 nT + M^2 |x_2|^2 + \frac{M^2 |\sigma_2|^2}{2\gamma} + K(T)\int_0^t \left( 1 + \frac{1}{\epsilon} e^{-\frac{2\gamma(t-s)}{\epsilon}}\right) \E \left[ |\overline{X}^\epsilon_1(s)|^2 + |\overline{X}^\epsilon_2(s)|^2 \right]ds.
\end{split}
\end{equation}
Hence, by applying Gronwall's inequality (see Theorem 15 of \cite{dragomir2003some}), we have that for all $t\in [0,T]$
\begin{align*}
&\E \left[|\overline{X}^\epsilon_1(t)|^2  +  |\overline{X}^\epsilon_2(t)|^2 \right]\\
&\leq \left(4|x_1|^2 + 4|\sigma_1|^2 nT + M^2 |x_2|^2 + \frac{M^2 |\sigma_2|^2}{2\gamma}\right) \exp\left[ K(T)\left(t + \frac{1}{2\gamma}\left(1 - e^{-\frac{2\gamma t}{\epsilon}}\right)\right)\right].
\end{align*}
Thus, for fixed $\epsilon \in (0,1]$
\begin{align*}
\sup_{t\in [0,T]}\E \left[ |\overline{X}^\epsilon_1(t)|^2  +  |\overline{X}^\epsilon_2(t)|^2 \right] &\leq \left(4|x_1|^2 + 4|\sigma_1|^2 nT + M^2 |x_2|^2 + \frac{M^2 |\sigma_2|^2}{2\gamma}\right) \exp\left[ K(T)\left(T + \frac{1}{2\gamma}\right)\right]\\
&:= K'(T).
\end{align*}
Since $K'(T)$ is independent of $\epsilon$, we obtain
\[
\sup_{\epsilon \in (0,1]} \sup_{t\in [0,T]}\E \left[ |\overline{X}^\epsilon_1(t)|^2  +  |\overline{X}^\epsilon_2(t)|^2 \right] \leq K'(T).
\]
\qed

The following theorem demonstrates that the objective function at the approximate feedback saddle point \eqref{approxoptimcontrol} is near the value of Problem (SLQG) with order $O(\epsilon)$. 

\begin{theorem}\label{theorem: approxvalue}
Suppose that Assumptions \ref{assumption: tech}, \ref{assumption: main}, \ref{assumption: negativedefinite} and \ref{assumption: attraction} hold. Let $T > 0$ be any finite time horizon and $\epsilon^*$ be the small positive parameter defined in Theorem \ref{tikhonov}. Then for all $0 < \epsilon <\epsilon^*$
\begin{equation}
J^\epsilon (x_1,x_2; \widehat{u}^\epsilon_1,\widehat{u}^\epsilon_2) - J^\epsilon (x_1,x_2; \overline{u}^\epsilon_1,\overline{u}^\epsilon_2) = O(\epsilon).
\end{equation}
\end{theorem}
\noindent \textit{Proof.} 
Similar to the proof in Theorem \ref{theorem: solution}, we apply Ito's lemma to $\left\langle P^\epsilon(t) X^\epsilon(t), X^\epsilon(t) \right\rangle$ where $X^\epsilon$ is defined in \eqref{statecompact} and by a completion of squares, we have that
\begin{align*}
&J^\epsilon (x_1,x_2; u_1,u_2)\\
&= - \frac{1}{2}\E\int_0^T\left[ |u_1(t) - \widehat{F}^\epsilon_{11}(t) X_1(t) - \widehat{F}^\epsilon_{12}(t) X^\epsilon_2(t) |^2 + |u_2(t) - \widehat{F}^\epsilon_{21}(t) X_1(t) - \widehat{F}^\epsilon_{22}(t) X^\epsilon_2(t) |^2 \right]dt\\
& \quad + \frac{1}{2}\left\langle P^\epsilon(0)x,x\right\rangle  + \frac{1}{2}\int_0^T\left\langle P^\epsilon(t) \sigma^\epsilon, \sigma^\epsilon\right\rangle dt.
\end{align*}
Thus, applying the approximate feedback saddle point \eqref{approxoptimcontrol}, we obtain
\begin{align*}
&|J^\epsilon (x_1,x_2; \overline{u}^\epsilon_1,\overline{u}^\epsilon_2) - J^\epsilon (x_1,x_2; \widehat{u}^\epsilon_1,\widehat{u}^\epsilon_2)|\\
&\leq \E\int_0^T\left[ \Big|(\overline{F}_{11}(t) - \widehat{F}^\epsilon_{11}(t)) \overline{X}^\epsilon_1(t)\Big|^2 + \Big|(\overline{F}_{12}(t) - \widehat{F}^\epsilon_{12}(t)) \overline{X}^\epsilon_2(t) \Big|^2  \right]dt\\
&\quad + \E \int_0^T \left[ \Big|(\overline{F}_{21}(t) - \widehat{F}^\epsilon_{21}(t)) \overline{X}^\epsilon_1(t)\Big|^2 +\Big| (\overline{F}_{22} - \widehat{F}^\epsilon_{22}(t)) \overline{X}^\epsilon_2(t) \Big|^2\right] dt.
\end{align*}
The result then follows from an application of Lemma \ref{lemma: feedbackopertorconverge} and Lemma \ref{lemma: boundedtraj}.
\qed 

The follow theorem gives an expression for the limiting value function of Problem (SLQG). 
\begin{theorem}\label{theorem: limitingvalue}
Suppose that Assumptions \ref{assumption: tech}, \ref{assumption: main}, \ref{assumption: negativedefinite} and \ref{assumption: attraction} hold. Let $T > 0$ be any finite time horizon and $\epsilon^*$ be the small positive parameter defined in Theorem \ref{tikhonov}. Define
\[
\overline{V}(x_1,x_2) = \frac{1}{2} \left\langle \overline{P}_{11}(0)x_1,x_1\right\rangle + \frac{1}{2} \int_0^T \langle \overline{P}_{11}(t) \sigma_1 , \sigma_1 \rangle dt + \frac{T}{2} \langle \overline{P}_{22} \sigma_2,  \sigma_2\rangle.
\]
Then for all $0<\epsilon < \epsilon^*$,
\begin{equation}
J^\epsilon (x_1,x_2; \widehat{u}^\epsilon_1,\widehat{u}^\epsilon_2) - \overline{V}(x_1,x_2) = O(\epsilon).
\end{equation}
\end{theorem}
\noindent \textit{Proof.} From Theorem \ref{theorem: solution},
\begin{align*}
J^\epsilon (x_1,x_2; \widehat{u}^\epsilon_1,\widehat{u}^\epsilon_2) - \overline{V}(x_1,x_2) &= \frac{1}{2}\left\langle P^\epsilon_{11}(0)x_1,x_1\right\rangle + \frac{\epsilon}{2} \left[ 2 \left\langle x_1, P^\epsilon_{12}(0) x_2\right\rangle + \left\langle P^\epsilon_{22}(0) x_2,x_2 \right\rangle \right]\\
&\quad + \frac{1}{2} \int_0^T \left[ \langle P_{11}^\epsilon(t) \sigma_1,  \sigma_1\rangle + \langle \overline{P}_{22}^\epsilon(t) \sigma_2,  \sigma_2\rangle \right] dt\\
&\quad - \frac{1}{2}\left\langle \overline{P}_{11}(0) x,x \right\rangle - \frac{1}{2} \int_0^T \langle \overline{P}_{11}(t) \sigma_1,  \sigma_1\rangle dt - \frac{T}{2} \langle \overline{P}_{22} \sigma_2,\sigma_2\rangle.
\end{align*}
From Theorem \ref{tikhonov},  we can simplify the above as
\begin{align*}
|J^\epsilon (x_1,x_2; \widehat{u}^\epsilon_1,\widehat{u}^\epsilon_2) - \overline{V}(x_1,x_2)| &\leq \epsilon K(x_1,x_2,T) + \frac{1}{2}\int_0^T \Big|\left\langle (\overline{P}_{22}(t) - \overline{P}_{22})\sigma_2,  \sigma_2 \right\rangle  \Big| dt
\end{align*}
where $K(x_1,x_2,T)$ is a constant which depends on the parameters $x_1,x_2,T$. By the Cauchy-Schwartz inequality, we have that
\begin{align*}
|J^\epsilon (x_1,x_2; \widehat{u}^\epsilon_1,\widehat{u}^\epsilon_2) - \overline{V}(x_1,x_2)| &\leq \epsilon K(x_1,x_2,T) + \frac{|\sigma_2|^2}{2}\int_0^T |\overline{P}_{22}(t) - \overline{P}_{22}| dt.
\end{align*}
Finally, an application of Corollary \ref{corollary: discbounds} gives the desired result.
\qed

\begin{remark}
It should be pointed out that two other common formulations of differential games are the deterministic case and also the stochastic case with multiplicative but no additive noise, see for example \cite{bernhard1979linear,mcasey2006generalized,yu2015optimal}. In both these cases, the value of the differential game is given as $\frac{1}{2} \langle P^\epsilon(0)x,x\rangle$ and subsequently, the limit of the value is $\frac{1}{2}\langle \overline{P}_{11}(0) x_1,x_1 \rangle$, which is devoid of the fast component of the differential game. As we see in the above result, when additive noise is present, the fast component $\overline{P}_{22}$ of the differential game makes a contribution to the limiting value.
\end{remark}

\section{One dimensional example}\label{section: onedimension}
In this section, we highlight our results for Problem (SLQG) when $k_i = m_i = n_i = 1$ for $i = 1,2$. The main observation in this case is that the solvability of the reduced differential and algebraic Riccati equations \eqref{reducedDRE}-\eqref{reducedARE} 
\begin{subnumcases}{}
\frac{d\overline{P}_{11}}{dt} + 2 \widetilde{A} \overline{P}_{11} + M \overline{P}_{11}^2 + N = 0,\quad \overline{P}_{11}(T) = 0, \label{reducedDRE1d} \\
\Delta_2 \overline{P}_{22}^2 + 2A_{22} \overline{P}_{22} + Q_2 = 0,\label{reducedARE1d}
\end{subnumcases}
where
\begin{align*}
\begin{cases}
\widetilde{A} = A_{11} + \frac{ \Delta Q_2 A_{21}  - A_{12} A_{21}  A_{22} }{A_{22}^2 - \Delta_2 Q_2}, \\
M = \Delta_1 + \frac{\Delta_2 A_{12}^2  - 2\Delta A_{12} A_{22} + \Delta^2 Q_2}{A_{22}^2 - \Delta_2 Q_2}, \\
N = Q_1 + \frac{ Q_2  A_{21}^2}{ A_{22}^2 - \Delta_2 Q_2},
\end{cases}
\end{align*}
has been well studied. Proposition 6.6.1 of \cite{jiongmin2014differential} states that reduced differential Riccati equation \eqref{reducedDRE1d} admits an unique solution $\overline{P}_{11} \in C([0,T];\R)$ if and only if either condition hold:
\begin{equation}\label{1dcondition1}
\widetilde{A}^2 -  M N \geq 0, \quad \text{ or }\quad  \widetilde{A} - \sqrt{| \widetilde{A}^2 -  M N |} \leq 0.
\end{equation}
For the reduced algebraic equation \eqref{reducedARE1d}, it is sufficient to assume that 
\begin{equation}\label{1dcondition2}
\Delta_2 Q_2 < 0.
\end{equation} 
In this case, it is straightforward to check that there exists an unique stabilising solution 
\[
\overline{P}_{22} = \frac{-A_{22} - \sqrt{A_{22}^2 - \Delta_2 Q_2}}{\Delta_2}
\]
and an unstable solution 
\[
\overline{P}_{22}^u = \frac{-A_{22} + \sqrt{A_{22}^2 - \Delta_2 Q_2}}{\Delta_2}.
\]
Under \eqref{1dcondition1} and \eqref{1dcondition2}, the requirements of Assumptions \ref{assumption: tech}, \ref{assumption: main} and \ref{assumption: negativedefinite} are satisfied. To see that Assumption \ref{assumption: attraction} is satisfied, let us begin by setting $\gamma := -S$. Then
\begin{align*}
\mathcal{R}_\delta = 
\begin{cases}
(-\infty, \frac{\delta - S}{\Delta_2}],\quad &\text{ if } \Delta_2 > 0,\\
[ \frac{\delta - S}{\Delta_2},\infty),\quad &\text{ if } \Delta_2 < 0.
\end{cases}
\end{align*}
We can lengthen the interval by heuristically considering the limit $\delta \rightarrow \gamma$ of $\mathcal{R}_\delta$ as
\begin{align*}
\lim_{\delta \rightarrow \gamma} \mathcal{R}_\delta=
\begin{cases}
(-\infty, \frac{-2S }{\Delta_2}),\quad &\text{ if } \Delta_2 > 0,\\
( \frac{-2S }{\Delta_2},\infty),\quad &\text{ if } \Delta_2 < 0.
\end{cases}
\end{align*}
Note here that the stable equilibrium $0$ of the boundary-layer problem is contained in $\mathcal{R}_\delta$ and the unstable equilibrium $\overline{P}_{22}^u - \overline{P}_{22} = \frac{-2S}{\Delta_2}$ of the boundary-layer problem is the boundary of $\lim_{\delta \rightarrow \gamma} \mathcal{R}_\delta$. This implies that $\lim_{\delta \rightarrow \gamma} \mathcal{R}_\delta$ is the exact region of attraction for the boundary-layer problem and furthermore, contains the values $\overline{P}_{22}$ and $-\overline{P}_{22}$. Hence there exists some $\delta \in (0,\gamma)$ such that $\overline{P}_{22}$ and $-\overline{P}_{22}$ are contained $\mathcal{R}_\delta$ \--- satisfying Assumption \ref{assumption: attraction}.

Thus, by Theorem \ref{theorem: maincorollary}, for any $T > 0$, there exists $\epsilon^* \in (0,1]$ such that for all $\epsilon \in (0,\epsilon^*)$ the generalised Riccati equation 
\begin{align*}
\begin{cases}
\dot{P}^\epsilon + (A^\epsilon)^* P^\epsilon + P^\epsilon A^\epsilon - P^\epsilon B^\epsilon R^{-1} (B^{\epsilon})^* P^\epsilon  + Q= 0,\\
P^\epsilon(T) = 0.
\end{cases}
\end{align*}
admits a unique solution $P^\epsilon \in C([0,T];\mathbb{S}^2)$. Consequently, the existence and uniqueness of a feedback saddle point and value follow from Theorem \ref{theorem: solution}, which also possess the asymptotic properties from Section \ref{section: approxcontrol}.

\section{Proof of Theorem \ref{tikhonov}}\label{section: proof}

For convenience, we begin by applying the change of variable $Z^\epsilon_1(t) = P_{12}^\epsilon(t) - h(P_{11}^\epsilon(t))$ and $Z^\epsilon_2(t) = P_{22}^\epsilon(t) - \overline{P}_{22}$ to the full system \eqref{full} and reformulating it as a initial value problem
\begin{subnumcases}{\label{shiftedfullsystem}}
\frac{dP_{11}^\epsilon}{dt} = f(P^\epsilon_{11},Z^\epsilon_1 + h(P^\epsilon_{11}), Z^\epsilon_2 + \overline{P}_{22}, \epsilon), \quad P_{11}^\epsilon(0) = 0,\\
\begin{split}
\epsilon\frac{d Z^\epsilon_1}{dt} &= g_1(P^\epsilon_{11},Z^\epsilon_1 + h(P^\epsilon_{11}), Z^\epsilon_2 + \overline{P}_{22}, \epsilon) \\
&\quad - \epsilon \frac{\partial h}{\partial P_{11}}(P^\epsilon_{11}) f(P^\epsilon_{11},Z^\epsilon_1 + h(P^\epsilon_{11}), Z^\epsilon_2 + \overline{P}_{22}, \epsilon), \quad Z^\epsilon_1(0) = -h(0),
\end{split}\\
\epsilon \frac{dZ^\epsilon_2}{dt} =g_2(P^\epsilon_{11},Z^\epsilon_1 + h(P^\epsilon_{11}), Z^\epsilon_2 + \overline{P}_{22}, \epsilon), \quad Z^\epsilon_2(0) = -\overline{P}_{22}.
\end{subnumcases}
Here we have used a slight abuse of notation in the time variable $t$ to change a terminal value problem to an initial value problem. i.e. $t = T-t$. Recall, from Assumption \ref{assumption: main} and Theorem \ref{theorem: reducedRiccati} that $\overline{P}_{11}(t)$ is the unique solution of
\[
\frac{d\overline{P}_{11}}{dt} = f(\overline{P}_{11},h(\overline{P}_{11}),\overline{P}_{22}, 0), \quad \overline{P}_{11}(0) = 0.
\]

Let $p_0$ be a positive constant such that $|\overline{P}_{11}(t)| \leq p_0$ for all $t\in [0,T]$ and define the closed ball $B_p := \{ P_{11} \in \mathbb{S}^{n_1}\ |\ |P_{11}(t)| \leq p,\ \forall t\in [0,T]\}$ where $p$ is arbitrarily chosen such that $p > p_0$. In order to show the trajectory $Z^\epsilon_1(t)$ does not blow up on the interval $[0,T]$, we will  require that $P_{11}^\epsilon$ is contained in the closed ball $B_p$. However, this compactness result has yet to be proven. To work around this, we  consider the smooth function $\psi: \mathbb{S}^{n_1} \rightarrow [0,1]$ defined such that $\psi(P_{11}) = 1$ when $|P_{11}| \leq \frac{1}{2}(p_0 + p)$ and $\psi(P_{11}) = 0$ when $|P_{11}| \geq p$. By heuristically replacing $P_{11}$ with $P_{11} \psi(P_{11})$ in the full system \eqref{shiftedfullsystem}, we obtain the following modified full system
\begin{subnumcases}{\label{augmentedfullsystem}}
\frac{dP_{11}^\epsilon}{dt} = \widetilde{f}(P_{11}^\epsilon, Z^\epsilon_1, Z^\epsilon_2, \epsilon), \quad P_{11}^\epsilon(0) = 0,\\
\epsilon \frac{dZ^\epsilon_1}{dt} = \widetilde{g}_1(P_{11}^\epsilon, Z^\epsilon_1, Z^\epsilon_2, \epsilon), \quad Z^\epsilon_1(0) = -h(0),\\
\epsilon \frac{dZ^\epsilon_2}{dt} = \widetilde{g}_2(P_{11}^\epsilon, Z^\epsilon_1, Z^\epsilon_2, \epsilon), \quad Z^\epsilon_2(0) = -\overline{P}_{22},
\end{subnumcases}
where the functions $\widetilde{f}, \widetilde{g}_1$ and $\widetilde{g}_2$ are defined as
\begin{align*}
\widetilde{f}(P_{11},Z_1,Z_2,\epsilon)
&= f(P_{11} \psi(P_{11}),Z_1 + h(P_{11} \psi(P_{11})), Z_2 + \overline{P}_{22}, \epsilon)\\
\widetilde{g}_1(P_{11},Z_1,Z_2,\epsilon) 
&= g_1(P_{11} \psi(P_{11}),Z_1 + h(P_{11} \psi(P_{11})), Z_2 + \overline{P}_{22}, \epsilon) \\
&\quad - \epsilon \frac{\partial h}{\partial P_{11}}(P_{11} \psi(P_{11})) f(P_{11} \psi(P_{11}), Z_1 + h(P_{11} \psi(P_{11})), Z_2 + \overline{P}_{22}, \epsilon)\\
\widetilde{g}_2(P_{11},Z_1,Z_2,\epsilon) 
&= g_2(P_{11} \psi(P_{11}),Z_1 + h(P_{11} \psi(P_{11})), Z_2 + \overline{P}_{22}, \epsilon).
\end{align*}
It is clear that when $|P_{11}| \leq \frac{1}{2} (p_0 + p)$, the modified full system \eqref{augmentedfullsystem} is identical to the original full system \eqref{shiftedfullsystem}. Moreover, as $\overline{P}_{11}(t)$ satisfies $|\overline{P}_{11}(t)| < \frac{1}{2} (p_0 + p)$ for all $t\in [0,T]$, it is also the solution of
\[
\frac{d\overline{P}_{11}}{dt} = \widetilde{f}(\overline{P}_{11},0,0,0),\quad \overline{P}_{11}(0) = 0.
\]
It is apparent that the function $\psi$ is chosen to ensure that $|P_{11} \psi(P_{11})|$ is bounded by $p$ for all $P_{11} \in \mathbb{S}^{n_1}$. Thus, if we find a solution to the modified full system $\eqref{augmentedfullsystem}$ with $|P_{11}^\epsilon(t)| \leq \frac{1}{2}(p_0 + p)$ for all $t\in [0,T]$, then it is also a solution to the original full system \eqref{shiftedfullsystem}. Moreover, by the locally Lipschitz property of $f,g_1$ and $g_2$ and the boundedness of $\partial h/\partial P_{11}$ on $B_p$, the solution to \eqref{augmentedfullsystem} is also the unique solution to the original full system \eqref{shiftedfullsystem}.

The remainder of the proof can be outlined in four parts: For sufficiently small $\epsilon,$
\begin{enumerate}[label=\Roman*.]
\item The continuous solutions $(Z^\epsilon_1, Z^\epsilon_2)$ are contained within the compact region $B_{2q_1} \times B_{\delta, q_2}$, for all $t\in [0,T]$, uniformly in $P_{11} \in \mathbb{S}^{n_1}$. In other words, $(Z^\epsilon_1, Z^\epsilon_2)$ is well-defined on $[0,T]$. Here, we note that $q_1$ can be chosen to be arbitrarily large.
\item The solution $P_{11}^\epsilon(t)$ is well-defined and continuous on $[0,T]$, and converges to $\overline{P}_{11}(t)$, for all $t\in [0,T]$, at a rate of $O(\epsilon)$. Moreover, $|P^\epsilon_{11}(t)| \leq \frac{1}{2} (p_0 + p)$ holds for all $t\in[0,T]$. Note that once this part is proven, the original and modified full systems will have the same solutions.
\item For all $t\in [0,T]$, the solutions $(Z^\epsilon_1(t), Z^\epsilon_2(t))$ converge to the solution $(\widehat{P}_{12}(t/\epsilon),\widehat{P}_{22}(t/\epsilon))$ of the boundary-layer problem \eqref{boundarylayer} at a rate $O(\epsilon)$.
\item Using the Lipschitz property of $h$, we can replace $h(P_{11}^\epsilon)$ with $h(\overline{P}_{11})$ to complete the proof.
\end{enumerate}

\textbf{Part I: Existence and uniqueness of $Z^\epsilon_1,Z^\epsilon_2$} 
 
In this part, we will show the existence and uniqueness of $Z^\epsilon_1$ and $Z^\epsilon_2$ in \eqref{augmentedfullsystem}. That is, the solution to
\begin{equation}\label{part2system}
\begin{split}
\frac{dZ^\epsilon_1}{dt} &= \frac{1}{\epsilon} \widetilde{g}_1(P_{11}^\epsilon,Z^\epsilon_1, Z^\epsilon_2, \epsilon), \quad Z^\epsilon_1(0) = -h(0),\\
\frac{dZ^\epsilon_2}{dt} &= \frac{1}{\epsilon} \widetilde{g}_2(P_{11}^\epsilon, Z^\epsilon_1, Z^\epsilon_2, \epsilon), \quad Z^\epsilon_2(0) = -\overline{P}_{22}.
\end{split}
\end{equation}
Adding and subtracting $\widetilde{g}_i(P_{11}^\epsilon, Z^\epsilon_1, Z^\epsilon_2, 0)$, we have that
\begin{align*}
\frac{dZ^\epsilon_1}{dt} &= \frac{1}{\epsilon} \widetilde{g}_1(P_{11}^\epsilon, Z^\epsilon_1, Z^\epsilon_2, 0) + \frac{1}{\epsilon} \left[\widetilde{g}_1(P_{11}^\epsilon, Z^\epsilon_1, Z^\epsilon_2, \epsilon) - \widetilde{g}_1(P_{11}^\epsilon, Z^\epsilon_1, Z^\epsilon_2, 0) \right], \quad Z^\epsilon_1(0) = -h(0),\\
\frac{dZ^\epsilon_2}{dt} &= \frac{1}{\epsilon} \widetilde{g}_2(P_{11}^\epsilon, Z^\epsilon_1, Z^\epsilon_2, 0) + \frac{1}{\epsilon} \left[\widetilde{g}_2(P_{11}^\epsilon, Z^\epsilon_1, Z^\epsilon_2, \epsilon) - \widetilde{g}_2(P_{11}^\epsilon, Z^\epsilon_1, Z^\epsilon_2, 0) \right], \quad Z^\epsilon_2(0)= -\overline{P}_{22}.
\end{align*}
For all $(Z^\epsilon_1, Z^\epsilon_2) \in B_{2q_1} \times B_{\delta,q_2}$, the functions $\widetilde{g}_i$ are Lipschitz with respect to $\epsilon$. That is, exists positive constants $M_1,M_2$, which may depend on $p,q_1,q_2,\delta$ and $T$, such that
\begin{align*}
|\widetilde{g}_1(P_{11}^\epsilon, Z^\epsilon_1, Z^\epsilon_2, \epsilon) - \widetilde{g}_1(P_{11}^\epsilon, Z^\epsilon_1, Z^\epsilon_2, 0)| \leq \epsilon M_1,\\
|\widetilde{g}_2(P_{11}^\epsilon, Z^\epsilon_1, Z^\epsilon_2, \epsilon) - \widetilde{g}_2(P_{11}^\epsilon, Z^\epsilon_1, Z^\epsilon_2, 0)| \leq \epsilon M_2.
\end{align*}
Then we have that
\begin{align*}
\frac{d|Z^\epsilon_1|^2}{dt} &\leq  \frac{2}{\epsilon}  \left\langle Z^\epsilon_1, Z^\epsilon_1 S + \Phi(P_{11}^\epsilon \psi(P_{11}^\epsilon)) Z^\epsilon_2 + Z^\epsilon_1 \Delta_2 Z^\epsilon_2\right\rangle +2 M_1 |Z^\epsilon_1|,\\
\frac{d|Z^\epsilon_2|^2}{dt} &\leq \frac{2}{\epsilon} \left\langle Z^\epsilon_2, S^* Z^\epsilon_2 + Z^\epsilon_2 S + Z^\epsilon_2 \Delta_2 Z^\epsilon_2\right\rangle + 2 M_2 |Z^\epsilon_2|,
\end{align*}
where $\Phi: \mathbb{S}^{n_1} \rightarrow \R^{n_1 \times n_2}$ is the mapping defined as
\begin{equation}
\Phi(P_{11}) = A_{21}^* + P_{11} \Delta + h(P_{11}) \Delta_2.
\end{equation}
Note that the function $\Phi(P_{11} \psi(P_{11}))$ is uniformly bounded with respect to $P_{11} \in \mathbb{S}^{n_1}$, say by $c_p$. For all $(Z^\epsilon_1, Z^\epsilon_2) \in B_{2q_1} \times B_{\delta,q_2}$, we can apply the Cauchy-Schwartz inequality to obtain
\begin{align*}
\frac{d|Z^\epsilon_1|^2}{dt} &\leq - \frac{2}{\epsilon}\left( \gamma - |\Delta_2| |Z^\epsilon_2| \right) |Z^\epsilon_1|^2 + \frac{2c_p }{\epsilon}  |Z^\epsilon_1| |Z^\epsilon_2| + 2M_1| Z^\epsilon_1|, \\
\frac{d|Z^\epsilon_2|^2}{dt} &\leq -\frac{2(\gamma - \delta)}{\epsilon} |Z^\epsilon_2|^2  + 2 M_2 |Z^\epsilon_2|.
\end{align*}
Applying the chain rule, we obtain
\begin{equation}\label{temps1}
\frac{d|Z^\epsilon_1|}{dt} \leq - \frac{1}{\epsilon}\left( \gamma - |\Delta_2| |Z^\epsilon_2| \right) |Z^\epsilon_1| + \frac{c_p}{\epsilon}   |Z^\epsilon_2| + M_1,
\end{equation}
\begin{equation}\label{temps2}
\frac{d|Z^\epsilon_2|}{dt} \leq -\frac{(\gamma - \delta)}{\epsilon} |Z^\epsilon_2|  +M_2.
\end{equation}
Hence \eqref{temps2} implies that
\begin{equation}\label{temps3a}
|Z^\epsilon_2(t)| \leq  e^{-\frac{(\gamma - \delta)t}{\epsilon} } |Z^\epsilon_2(0)| + \frac{\epsilon M_2}{\gamma - \delta} \left(1 - e^{-\frac{(\gamma -\delta) t}{\epsilon}} \right), \quad \forall t\in [0,T]
\end{equation}
and simplifies to
\begin{equation}\label{temps3}
|Z^\epsilon_2(t)| \leq  e^{-\frac{(\gamma - \delta)t}{\epsilon} } |Z^\epsilon_2(0)| + \frac{\epsilon M_2}{\gamma - \delta}, \quad \forall t\in [0,T].
\end{equation}
By Assumption \ref{assumption: attraction}, $Z^\epsilon_2(0) = -\overline{P}_{22}$ is in the interior of the ball $B_{\delta,q_2}$ (i.e. $|\overline{P}_{22}| < q_2$). Hence, for sufficiently small $\epsilon$, the continuous solution $Z^\epsilon_2(t)$ is contained in the interior of the closed ball $B_{\delta,q_2}$ for all $t \in [0,T]$. On the other hand, substituting \eqref{temps3a} into \eqref{temps1} gives
\begin{align*}
\frac{d|Z^\epsilon_1|}{dt}
&\leq - \frac{1}{\epsilon}\left[ \gamma - |\Delta_2|q_2  e^{-\frac{(\gamma - \delta)t}{\epsilon} } - \frac{\epsilon M_2|\Delta_2|}{\gamma - \delta} \right] |Z^\epsilon_1|  + \frac{c_p}{\epsilon}  e^{-\frac{(\gamma - \delta)t}{\epsilon} } q_2 + M_1 + \frac{M_2 c_p}{\gamma - \delta} \left(1 -  e^{-\frac{(\gamma -\delta) t}{\epsilon}} \right)
\end{align*}
Hence
\begin{align*}
|Z^\epsilon_1(t)| &\leq  \exp\left( -\frac{\gamma t}{\epsilon} + \frac{M_2|\Delta_2| t}{\gamma - \delta} + \frac{|\Delta_2| q_2}{\gamma - \delta} \left(1 -  e^{-\frac{(\gamma - \delta)t}{\epsilon}} \right)\right) |Z^\epsilon_1(0)| \\
&\quad + \int_0^t  \exp\left( -\frac{\gamma (t - s)}{\epsilon} + \frac{M_2|\Delta_2| (t - s)}{\gamma - \delta} + \frac{|\Delta_2| q_2}{\gamma - \delta} \left( e^{-\frac{(\gamma - \delta)s}{\epsilon}} -  e^{-\frac{(\gamma - \delta)t}{\epsilon}} \right) \right) \\
&\qquad \left[\frac{c_p}{\epsilon}  e^{-\frac{(\gamma - \delta)s}{\epsilon} } q_2 + M_1 + \frac{M_2 c_p}{\gamma - \delta} \left(1 -  e^{-\frac{(\gamma -\delta) s}{\epsilon}} \right) \right] ds\\
&\leq e^{ \frac{|\Delta_2| q_2}{\gamma - \delta}} e^{\left(\frac{M_2|\Delta_2| }{\gamma - \delta} - \frac{\gamma}{\epsilon}\right) t} |Z^\epsilon_1(0)| \\
&\quad + e^{ \frac{|\Delta_2| q_2}{\gamma - \delta}} e^{\left(\frac{M_2|\Delta_2| }{\gamma - \delta} - \frac{\gamma}{\epsilon}\right) t} \int_0^t  e^{\frac{\gamma s}{\epsilon}} \left[\frac{c_p}{\epsilon}  e^{-\frac{(\gamma - \delta)s}{\epsilon} } q_2 + M_1 + \frac{M_2 c_p}{\gamma - \delta}  \right] ds\\
&= e^{ \frac{|\Delta_2| q_2}{\gamma - \delta}} e^{\left(\frac{M_2|\Delta_2| }{\gamma - \delta} - \frac{\gamma}{\epsilon}\right) t} |Z^\epsilon_1(0)| \\
&\quad + e^{ \frac{|\Delta_2| q_2}{\gamma - \delta}} e^{\left(\frac{M_2|\Delta_2| }{\gamma - \delta} - \frac{\gamma}{\epsilon}\right) t} \int_0^t  \left[\frac{c_p q_2}{\epsilon}  e^{-\frac{\delta s}{\epsilon} } + e^{\frac{\gamma s}{\epsilon}} \left( M_1 + \frac{M_2 c_p}{\gamma - \delta}  \right)  \right]ds. 
\end{align*}
Evaluating the integral gives
\begin{equation}\label{temps5}
|Z^\epsilon_1(t)| \leq e^{ \frac{|\Delta_2| q_2}{\gamma - \delta}} e^{\left(\frac{M_2|\Delta_2| }{\gamma - \delta} - \frac{\gamma}{\epsilon}\right) t} \left[|Z^\epsilon_1(0)| + \frac{c_p q_2}{\delta}\right]+ \epsilon e^{ \frac{|\Delta_2| q_2}{\gamma - \delta}} e^{\left(\frac{M_2|\Delta_2| }{\gamma - \delta} - \frac{(\gamma - \delta)}{\epsilon}\right) t} \left[ M_1 + \frac{M_2 c_p}{\gamma - \delta}  \right].
\end{equation}
Set $q_1 > e^{\frac{|\Delta_2| q_2}{\gamma - \delta}} \left(|h(0)| + \frac{c_p q_2}{\delta}\right) + |h(0)|$ and choose $\epsilon$ be sufficiently small such that
\[
\epsilon e^{ \frac{|\Delta_2| q_2}{\gamma - \delta}} \left[ M_1 + \frac{M_2 c_p}{\gamma - \delta}  \right] < q_1 \text{ and } \frac{M_2|\Delta_2| }{\gamma - \delta} - \frac{(\gamma - \delta)}{\epsilon} < 0.
\]
As a result
\begin{equation}
|Z^\epsilon_1(t)|
\leq e^{ \frac{|\Delta_2| q_2}{\gamma - \delta}} \left[|h(0)| + \frac{c_p q_2}{\delta}\right]+ \epsilon e^{ \frac{|\Delta_2| q_2}{\gamma - \delta}} \left[ M_1 + \frac{M_2 c_p}{\gamma - \delta}  \right] < 2q_1.
\end{equation}
Hence, for sufficiently small $\epsilon$, the continuous solution $Z^\epsilon_1(t)$ is contained in interior of the closed ball of radius $2q_1$ for all $t\in [0,T]$.

In this part, we have shown that, for sufficiently small $\epsilon$, if $(Z_1^\epsilon(t),Z_2^\epsilon(t)) \in B_{2q_1} \times B_{\delta,q_2}$ for all $t\in [0,T]$ then $(Z_1^\epsilon(t),Z_2^\epsilon(t))$ must lie in the interior of the closed set $B_{2q_1} \times B_{\delta,q_2}$ for all $t\in [0,T]$. As a result, for sufficiently small $\epsilon$, $(Z_1^\epsilon(t),Z_2^\epsilon(t))$ is indeed contained in the closed set $B_{2q_1} \times B_{\delta,q_2}$ for all $t \in [0,T]$. For otherwise, because of the continuity of $(Z_1^\epsilon,Z_2^\epsilon)$, there would exist $t\in (0,T]$ such that say $|Z_1^\epsilon(t)| = 2q_1$. However, as we have just established, this implies that $|Z_1^\epsilon(t)| < 2q_1$, which is a contradiction. \\

\textbf{Part II: Existence, uniqueness and convergence of $P^\epsilon_{11}$}

In this part, we denote $N_i$ as positive constants, which may depend on the parameters $p,q_1,q_2,\delta$ and $T$. Consider the slow component of the modified full system \eqref{augmentedfullsystem} given by
\begin{equation}
\frac{dP_{11}^\epsilon}{dt} = \widetilde{f}(P_{11}^\epsilon,Z^\epsilon_1,Z^\epsilon_2,\epsilon),\quad P_{11}^\epsilon(0) = 0.
\end{equation}
Recall that by Assumption \ref{assumption: main}, the function $\overline{P}_{11}(t), t\in [0,T]$ is a solution to the reduced equation
\begin{equation}\label{temp22}
\frac{d\overline{P}_{11}}{dt} = \widetilde{f}(\overline{P}_{11}, 0,0,0), \quad \overline{P}_{11}(0) = 0.
\end{equation}
Since, for all $(P_{11}, Z^\epsilon_1(t), Z^\epsilon_2(t)) \in \mathbb{S}^{n_1}\times B_{2q_1} \times B_{\delta, q_2}$, the function $\widetilde{f}$ is Lipschitz continuous with respect to its parameters, we have that
\begin{align*}
|\widetilde{f}(P_{11}^\epsilon,Z^\epsilon_1,Z^\epsilon_2,\epsilon) - \widetilde{f}(\overline{P}_{11}, 0,0,0)| &\leq \epsilon N_1 + N_2 | Z^\epsilon_1| + N_3 | Z^\epsilon_2| + N_4 | P_{11}^\epsilon - \overline{P}_{11}|.
\end{align*}
Let $y^\epsilon(t) = P_{11}^\epsilon(t) - \overline{P}_{11}(t)$. Then
\begin{align*}
y^\epsilon(t) &= \int_0^t \left[ \widetilde{f}(P_{11}^\epsilon(s),Z^\epsilon_1(s),Z^\epsilon_2(s),\epsilon) - \widetilde{f}(\overline{P}_{11}(s), 0,0,0) \right] ds.
\end{align*}
Taking the norm
\begin{align*}
|y^\epsilon(t)| &\leq \int_0^t | \widetilde{f}(P_{11}^\epsilon(s),Z^\epsilon_1(s),Z^\epsilon_2(s),\epsilon) - \widetilde{f}(\overline{P}_{11}(s), 0,0,0) | ds\\
&\leq \int_0^t \left[\epsilon N_1 + N_2 | Z^\epsilon_1(s)| + N_3 | Z^\epsilon_2(s)| + N_4 |y^\epsilon(s)| \right]ds.
\end{align*}
From \eqref{temps3} and \eqref{temps5}
\begin{align*}
|y^\epsilon(t)|
&\leq \int_0^t \left[\epsilon N_5 + N_6  e^{-\frac{(\gamma - \delta) s}{\epsilon}} + N_4 |y^\epsilon(s)| \right]ds\\
&= \epsilon N_5 t + \frac{\epsilon N_6 }{\gamma - \delta} \left( 1 - e^{-\frac{(\gamma - \delta) t}{\epsilon}}\right) + N_4 \int_0^t |y^\epsilon(s)| ds\\
&\leq \epsilon N_5 T + \frac{\epsilon N_6 }{\gamma - \delta}  + N_4 \int_0^t |y^\epsilon(s)| ds.
\end{align*}
Finally, Gronwall's inequality gives
\begin{equation}\label{tempw}
|P^\epsilon_{11}(t) - \overline{P}_{11}(t)| \leq \epsilon N_7 (1 + T) e^{N_4 T}
\end{equation}
Hence, for sufficiently small $\epsilon$, we obtain $|P^\epsilon_{11}(t)| < \frac{1}{2}(p_0 + p)$ for all $t\in [0,T]$. This implies that the continuous solution $P_{11}^\epsilon(t)$ exists for all $t\in [0,T]$, and satisfies the above estimate.\\

\textbf{Part III: Convergence of $Z^\epsilon_1,Z^\epsilon_2$}

In this part, $L_i$ denotes a positive constant, which may depend on the parameters $p,q_1,q_2,\delta$ and $T$. Let us work on the time scale $\tau = t/\epsilon$. Denote the difference $v^\epsilon_i(\tau) = Z_i^\epsilon(\epsilon\tau) - \widehat{P}_{i2}(\tau)$ for $i = 1,2$. Differentiating $v^\epsilon_2(\tau)$ with respect to $\tau$ gives
\begin{align*}
\frac{dv^\epsilon_2}{d\tau} &= \widetilde{g}_2(P^\epsilon_{11}, Z^\epsilon_1, Z^\epsilon_2, \epsilon) - \widetilde{g}_2(0, \widehat{P}_{12}, \widehat{P}_{22}, 0)\\
&= \widetilde{g}_2(P_{11}^\epsilon, Z^\epsilon_1, v_2^\epsilon, 0) + \Gamma_1 + \Gamma_2\\
&= S^* v_2^\epsilon + v_2^\epsilon S + v_2^\epsilon \Delta_2 v_2^\epsilon +  \Gamma_1 + \Gamma_2 
 \end{align*}
where
\begin{align*}
\Gamma_1 &= \widetilde{g}_2(P^\epsilon_{11}, Z^\epsilon_1, Z^\epsilon_2, 0) - \widetilde{g}_2(0, \widehat{P}_{12}, \widehat{P}_{22}, 0) - \widetilde{g}_2(P_{11}^\epsilon, Z^\epsilon_1, v_2^\epsilon, 0),\\
\Gamma_2 &= \widetilde{g}_2(P^\epsilon_{11}, Z^\epsilon_1, Z^\epsilon_2, \epsilon) - \widetilde{g}_2(P^\epsilon_{11}, Z^\epsilon_1, Z^\epsilon_2, 0).
\end{align*}
Note that the function $\widetilde{g}_2(P_{11}, Z_1, Z_2, 0)$ is independent of the variables $P_{11}$ and $Z_1$. By the Mean Value Theorem, we have that 
\[
\Gamma_1 = \left[ \frac{\partial \widetilde{g}_2}{\partial Z_2}(P^\epsilon_{11},Z^\epsilon_1, m_1 v_2^\epsilon + \widehat{P}_{22}, 0) - \frac{\partial \widetilde{g}_2}{\partial Z_2}(P^\epsilon_{11} ,Z^\epsilon_1, m_2 v_2^\epsilon, 0) \right] v_2^\epsilon
\]
where $0 <  m_1,m_2 < 1$. Since $ \partial \widetilde{g}_2/ \partial Z_2 $ is Lipschitz in $Z_2$, we have that
\[
|\Gamma_1| \leq L_2 |v_2^\epsilon|^2 + L_3 |v_2^\epsilon| |\widehat{P}_{22}|.
\]
Since $\widetilde{g}_2$ is Lipschitz with respect to $\epsilon$, we have that $|\Gamma_2| \leq \epsilon L_1$. Hence when $v_2^\epsilon \in B_{\delta,q_2}$ and $|v_2^\epsilon| \leq \frac{\gamma - \delta}{L_2}$
\begin{align*}
\frac{d|v_2^\epsilon|^2}{d\tau} &\leq -2(\gamma - \delta) |v_2^\epsilon|^2 + 2\epsilon L_1 |v_2^\epsilon| + 2L_2 |v_2^\epsilon|^3 + 2L_3 |v_2^\epsilon|^2 |\widehat{P}_{22}|\\
&\leq -(\gamma - \delta) |v_2^\epsilon|^2 + 2\epsilon L_1 |v_2^\epsilon| + 2L_3 q_2 |v_2^\epsilon|^2 e^{-(\gamma - \delta) \tau}\\
&= - \left( (\gamma - \delta) - 2L_3 q_2 e^{-(\gamma - \delta) \tau} \right) |v_2^\epsilon|^2 + 2\epsilon L_1 |v_2^\epsilon|.
\end{align*}
By the chain rule,
\begin{align*}
\frac{d|v_2^\epsilon|}{d\tau} \leq -  \left(\frac{1}{2} (\gamma - \delta) - L_3 q_2 e^{-(\gamma - \delta) \tau} \right) |v_2^\epsilon| + \epsilon L_1.
\end{align*}
Hence
\begin{align*}
|v_2^\epsilon(\tau)| &\leq \epsilon L_1  \int_0^\tau \exp \left[ -  \int^\tau_s \left( \frac{1}{2}(\gamma - \delta) - L_3 q_2 e^{-(\gamma - \delta) r} \right)dr \right] ds\\
&\leq \epsilon L_1 e^{\frac{L_3 q_2}{\gamma - \delta}} \int_0^\tau e^{-\frac{1}{2} (\gamma - \delta)(\tau - s)} ds\\
&=  \frac{2\epsilon L_1}{\gamma - \delta} e^{\frac{L_3 q_2}{\gamma - \delta}} \left( 1 - e^{-\frac{1}{2} (\gamma - \delta) \tau} \right)
\end{align*}
which simplifies to
\begin{equation}\label{temp99}
|v_2^\epsilon(\tau)|  \leq \frac{2\epsilon L_1}{\gamma - \delta} e^{\frac{L_3 q_2}{\gamma - \delta}}.
\end{equation}
Changing the time variable back to $t = \epsilon \tau$, we have that
\begin{equation}\label{temp00}
|Z^\epsilon_2(t) - \widehat{P}_{22}(t/\epsilon)| \leq \frac{2\epsilon L_1}{\gamma - \delta} e^{\frac{L_3 q_2}{\gamma - \delta}}  ,\quad \forall t\in [0,T].
\end{equation}
Hence, when $\epsilon$ is chosen is small enough such that 
\begin{equation}\label{suffeps}
\frac{2\epsilon L_1}{\gamma - \delta} e^{\frac{L_3 q_2}{\gamma - \delta}} \leq \min\left\{q_2, \frac{\gamma - \delta}{L_2} \right\},
\end{equation}
the estimate \eqref{temp00} holds.

We can repeat a similar line of reasoning for $v^\epsilon_1(\tau)$. Differentiating $v^\epsilon_1(\tau)$ gives
\begin{align*}
\frac{dv^\epsilon_1}{d\tau} &= \widetilde{g}_1(P^\epsilon_{11}, Z^\epsilon_1, Z^\epsilon_2, \epsilon) - \widetilde{g}_1(0,\widehat{P}_{12}, \widehat{P}_{22}, 0)\\
&= \widetilde{g}_1(P_{11}^\epsilon, v_1^\epsilon, v_2^\epsilon, 0) + \Gamma_3 + \Gamma_4 + \Gamma_5\\
&= v_1^\epsilon S + \Phi(P^\epsilon_{11} \psi(P^\epsilon_{11})) v_2^\epsilon + v_1^\epsilon \Delta_2 v_2^\epsilon +  \Gamma_3 + \Gamma_4 + \Gamma_5
\end{align*}
where
\begin{align*}
\Gamma_3 &= \widetilde{g}_1(P^\epsilon_{11}, Z^\epsilon_1, Z^\epsilon_2, 0) - \widetilde{g}_1(P^\epsilon_{11}, \widehat{P}_{12}, \widehat{P}_{22}, 0) - \widetilde{g}_1(P^\epsilon_{11}, v_1^\epsilon, v_2^\epsilon, 0),\\
\Gamma_4 &= \widetilde{g}_1(P^\epsilon_{11}, \widehat{P}_{12}, \widehat{P}_{22}, 0) - \widetilde{g}_1(0, \widehat{P}_{12}, \widehat{P}_{22}, 0),\\
\Gamma_5 &= \widetilde{g}_1(P^\epsilon_{11}, Z^\epsilon_1, Z^\epsilon_2, \epsilon) - \widetilde{g}_1(P^\epsilon_{11}, Z^\epsilon_1, Z^\epsilon_2, 0).
\end{align*}
Here recall that $\widetilde{g}_1(P_{11},Z_1,Z_2,0) = Z_1 S + \Phi(P_{11} \psi(P_{11})) Z_2 + Z_1 \Delta Z_2$. Using this, we can show that
\[
|\Gamma_3| \leq |\Delta| \left( |v_1^\epsilon|| \widehat{P}_{22}| + |v_2^\epsilon| |\widehat{P}_{12}|\right)
\]
and since $P_{11}^\epsilon(t)$ is defined on $t\in [0,T]$,
\[
|\Gamma_4| \leq |\Phi(P^\epsilon_{11} \psi(P^\epsilon_{11})) - \Phi(0)| |\widehat{P}_{22}| \leq  L_4 |P_{11}^\epsilon(\epsilon \tau)| |\widehat{P}_{22}| \leq \epsilon \tau L_5 |\widehat{P}_{22}|.
\]
Lastly, $\widetilde{g}_1$ is Lipschitz with respect to $\epsilon$, which implies that $|\Gamma_5| \leq \epsilon L_6 $. Putting this together we can deduce that 
\begin{align*}
\frac{d|v^\epsilon_1|^2}{d\tau} 
&\leq -2\gamma |v_1^\epsilon|^2 + 2c_p |v_2^\epsilon| |v_1^\epsilon| + 2|\Delta| |v_2^\epsilon| |v_1^\epsilon|^2 \\
&\quad + 2|\Delta| \left( |v_1^\epsilon|| \widehat{P}_{22}| + |v_2^\epsilon| |\widehat{P}_{12}|\right) |v^\epsilon_1| + 2\epsilon \tau L_5 |\widehat{P}_{22}| |v_1^\epsilon|+ \epsilon L_6 |v_1^\epsilon|\\
&\leq 2 \left( |\Delta| |v_2^\epsilon| + |\Delta| |\widehat{P}_{22}|  - \gamma \right) |v_1^\epsilon|^2 + 2 (c_p + |\Delta||\widehat{P}_{12}|) |v_2^\epsilon| |v_1^\epsilon| + 2\epsilon \tau L_5 |\widehat{P}_{22}| |v_1^\epsilon|+ 2\epsilon L_6 |v_1^\epsilon|.
\end{align*}
By the chain rule, we have that
\begin{align*}
\frac{d|v^\epsilon_1|}{d\tau} \leq   \left( |\Delta| |v_2^\epsilon| + |\Delta| |\widehat{P}_{22}|  - \gamma \right) |v_1^\epsilon| +  (c_p + |\Delta||\widehat{P}_{12}|) |v_2^\epsilon|  + \epsilon \tau L_5 |\widehat{P}_{22}| + \epsilon L_6.
\end{align*}
From Lemma \ref{lemma: exponentialstabilising} and \eqref{temp99}, we have that
\begin{align*}
\frac{d|v^\epsilon_1|}{d\tau} &\leq   \left(\epsilon L_7 + |\Delta| e^{-(\gamma -\delta)\tau} q_2  - \gamma \right) |v_1^\epsilon| +  \epsilon L_8 + \epsilon \tau L_5 e^{-(\gamma - \delta) \tau} q_2\\
&\leq  \left(\epsilon L_7 + |\Delta| e^{-(\gamma -\delta)\tau} q_2  - \gamma \right) |v_1^\epsilon| +  \epsilon L_9.
\end{align*}
Here we used the fact that $\tau e^{-(\gamma-\delta)\tau} \leq (\gamma - \delta)^{-1}$. Hence, when the condition \eqref{suffeps} is satisfied,
\begin{equation}\label{tempv1}
\begin{split}
|v_1^\epsilon(\tau)| &\leq \epsilon L_9 \int_0^\tau \exp \left[ -  \int^\tau_s \left( \gamma -  \epsilon L_7 - |\Delta| q_2 e^{-(\gamma -\delta)r} \right) dr \right] ds\\
&\leq \frac{2\epsilon L_9}{\gamma} e^{\frac{|\Delta| q_2}{\gamma - \delta}}.
\end{split}
\end{equation}
\textbf{Part IV: Final estimates}\\
Using the Lipschitz property of $h$ and the inequalities \eqref{tempw} and \eqref{tempv1}, we have that
\begin{align*}
|P_{12}^\epsilon(t) - h(\overline{P}_{11}(t)) - \widehat{P}_{12}(t/\epsilon)| &= |v_1^\epsilon(t/\epsilon) + h(\overline{P}_{11}(t))- h(P^\epsilon_{11}(t))| \\
&\leq \frac{2\epsilon L_9}{\gamma} e^{\frac{|\Delta| q_2}{\gamma - \delta}} + L_{10} | \overline{P}_{11}^\epsilon(t) - \overline{P}_{11}(t)| \\
&\leq \epsilon  L_{11}
\end{align*}
where $L_{10}$ and $L_{11}$ are positive constants that may depend on the parameters $p,q_1,q_2,\delta$ and $T$.
\qed

\bibliographystyle{abbrv}
\bibliography{references}

\end{document}